\documentclass[12pt,leqno]{article}
\usepackage{amssymb}

\newcommand\qed{\hfill$\sqcap\kern-7.5pt\hbox{$\sqcup$}$}

\newcommand{\NN}{\mathbb{N}}
\newcommand{\RR}{\mathbb{R}}
\newcommand{\R}{\mathbb{R}}
\newcommand{\N}{\mathbb{N}}
\newcommand{\Sph}{\mathbb{S}}

\newtheorem{theo}{Theorem}
\newtheorem{prop}[theo]{Proposition}
\newtheorem{lem}[theo]{Lemma}
\newtheorem{cor}[theo]{Corollary}
\newtheorem{rem}[theo]{Remark}
\newtheorem{rems}[theo]{Remarks}
\newtheorem{defin}[theo]{Definition}
\renewcommand{\theequation}{\thesection.\arabic{equation}}
\renewcommand{\thetheo}{\thesection.\arabic{theo}}

\newcommand{\beqn}{\begin{equation}}
\newcommand{\eeqn}{\end{equation}}
\newcommand{\bear}{\begin{eqnarray}}
\newcommand{\eear}{\end{eqnarray}}
\newcommand{\bean}{\begin{eqnarray*}}
\newcommand{\eean}{\end{eqnarray*}}

\newcommand{\Cc}{\mathcal{C}}
\newcommand{\DD}{\mathcal{D}}
\newcommand{\EE}{\mathcal{E}}
\newcommand{\OO}{\mathcal{O}}

\newcommand{\e}{{\varepsilon}}

\newcommand{\eps}{\varepsilon}
\newcommand{\wto}{\rightharpoonup}

\begin{document}

\title{Cooling process for inelastic Boltzmann equations for hard
spheres, Part I: The Cauchy problem}

\author{ S. {\sc Mischler}$^1$, C. {\sc Mouhot}$^2$
and M. {\sc Rodriguez Ricard}$^3$}

\footnotetext[1]{Ceremade, Universit\'e Paris IX-Dauphine,
place du M$^{al}$ DeLattre de Tasigny, 75016 Paris, France.}

\footnotetext[2]{UMPA, \'ENS Lyon, 46, alle d'Italie
69364 Lyon cedex 07, France.}

\footnotetext[3]{Facultad de Matem\'atica y Computaci\'on, Universidad
de La Habana, C.Habana 10400, Cuba.}

\maketitle

\begin{abstract}
We develop the Cauchy theory of the spatially
homogeneous inelastic Boltzmann equation for hard spheres, for
a general form of collision rate which includes in particular
variable restitution coefficients depending on the kinetic energy and
the relative velocity as well as the sticky particles model.
We prove (local in time) non-concentration estimates in Orlicz spaces,
from which we deduce weak stability and existence theorem.
Strong stability together with uniqueness and instantaneous appearance
of exponential moments are proved  under
additional smoothness assumption on the initial datum,
for a restricted class of collision rates.
Concerning the long-time behaviour, we give conditions for
the cooling process to occur or not in finite time.
\end{abstract}



\tableofcontents

\vspace{0.3cm}

\section{Introduction and main results}

\setcounter{equation}{0}
\setcounter{theo}{0}

In this paper we address the Cauchy problem for the spatially
homogeneous Boltzmann equation modelling the dynamic of a homogeneous
system of
inelastic hard spheres which interact only through binary collisions.
More precisely, describing the gas by the probability density
$f(t,v) \ge 0$ of particles with velocity $v \in \RR^N$ ($N \ge 2$)
at time $t \ge 0$, we study the existence, uniqueness and the
qualitative behaviour
of solutions to the Boltzmann equation for inelastic collision
     \bear  \label{eqB1}
     {\partial f \over \partial t}& = & Q(f,f) \quad\hbox{ in }\quad
(0,+\infty)\times \RR^N,\\
     \label{eqB2}
     f(0,\cdot) & = &  f_{\mbox{\scriptsize{in}}} \quad\hbox{ in }\quad
\RR^N.
     \eear
The use of Boltzmann inelastic hard spheres-like models to describe
dilute, rapid flows of granular media started with the seminal physics
paper \cite{Haff}, and a huge physics litterature has
developed in the last twenty years. The study of granular systems in
such regime is motivated by
their unexpected physical behavior (with the phenomena of {\em
collapse} --or ``{\em cooling effect}''--
at the kinetic level and {\em clustering} at the hydrodynamical level),
their use to derive hydrodynamical
equations for granular fluids, and their applications.

\smallskip

From the mathematical viewpoint, works on the Cauchy problem for
Boltzmann models have been first restricted to the so-called {\em
inelastic Maxwell molecules model}, where, in a similar way to the 
Maxwell model in the elastic framework, the collision rate is independent on the
relative velocity. Existence, uniqueness of solutions and description
of the asymptotic cooling has been obtained in  \cite{BCG00}  for the
{\em inelastic Maxwell molecules model} with constant normal
restitution coefficients as well as with
some cases of normal restitution coefficients depending on the kinetic
energy of the solution. More precise properties of the solutions, 
such as their convergence to self-similarity, have also been 
investigated and we refer to our companion paper~\cite{MM05} 
for more details and more references on this issue.  
At least in the spatially homogeneous setting, 
the {\em inelastic Maxwell molecules model} seems well understood now. 
The {\em Maxwell molecules model} is important because of its analytic
simplifications (with regards to the hard sphere model) allowing to use
powerful Fourier transform tools as introduced
by Bobylev (see for instance~\cite{BobyMaxwell}) for the {\em elastic Maxwell molecules
Boltzmann equation}. Another simplification which has lead to
interesting results is the restriction to one-dimensional models (in
space and velocity), where, on the contrary to the elastic case, the
collision operator has a non-trivial outcome. These models have been
considered in  \cite{BCP97,Toscani00,BP00} for some
cases of normal restitution coefficients possibly depending on the
relative velocity.

It is possible to modify the collision operator of the {\em
inelastic Maxwell molecules model} by a multiplication by a function of
the kinetic energy in order to restore its dimensional homogeneity (see
\cite{BCG00}) and thus the rate of cooling. Physically the derivation of 
this model amounts to replace the collision rate by a mean value 
independent on the relative velocity, starting from the 
{\em inelastic hard spheres model}, and the resulting approximation is named   
{\em pseudo-Maxwell molecules} in \cite{BCG00}. However, 
fine properties of the distribution (such as the behavior of the 
overpopulated tails or the self-similar solutions) are broken or modified
by that approximation with respect to the {\em inelastic hard spheres model}.
The recent papers \cite{GPV**,BGP**} have studied the case of inelastic
hard spheres with constant normal restitution coefficients in any
dimension and in various regimes: in particular in a thermal bath,
{\em i.e.}, when a heat source term is added to the equation, and in the
self-similar variables of the free regime.
Existence and smoothness of solutions to the Cauchy problem and to the
associated stationary problem are obtained in \cite{GPV**} for the
thermal bath regime, while precise estimates on the tails of the
stationary solutions (assuming their existence) 
for various regimes (including the two ones above-mentioned) 
are exhibited in \cite{BGP**}.

\smallskip

In the present work, we shall construct solutions to the freely cooling
Boltzmann equation for {\em inelastic hard spheres} in any dimension $N \ge 2$
and for a general framework of distributions of inelasticity (defined by a measure 
on the set of all possible post-collisional velocities), 
which covers in particular {\em variable normal restitution
coefficients} possibly depending on the relative velocity
and the kinetic energy of the solution. 
It includes the cases of {\em visco-elastic hard spheres model} (see \cite{BriPoe}) 
as well as the case of {\em sticky particles model}.
Our framework enables to consider interesting physical features, such
as elasticity increasing when the relative velocity or the temperature
decrease ({\em ``normal'' granular media}) or the opposite phenomenon
({\em ``anomalous'' granular media}). 
We refer to \cite{BCG00,Toscani00,FM**,BriPoe} and the references therein
for a physical motivation. Let us emphasize that these solutions are
new even in the case of a constant normal restitution  coefficient as
considered in \cite{GPV**,BGP**}. We also discuss the uniqueness of
solutions, the instantaneous appearance of exponential moments and various
conditions on the collisions rate for the collapse to occur or not in
finite time. A second part of this work \cite{MM05} will be concerned with the
existence of self-similar solutions
and the tail behavior of the distribution. In a third part
\cite{MM**}, we shall prove the uniqueness and the asymptotic stability
of these self-similar solutions for a small inelasticity.

\smallskip
From the viewpoint of {\it mathematical tools}, our main new
contributions can be summarized as follows:

\smallskip
(i) A generalization of the propagation of the $L^p$-norm of the
solution for the elastic Boltzmann equation based on Young's 
inequality as introduced in \cite{DM**} (see also \cite{MR03,ELM**}
where similar ideas are used for a different model), into a result of
propagation of Orlicz norms for inelastic (and elastic) Boltzmann equations. 
This {\em a priori} estimate is used in order to prove the existence of solutions to the
inelastic Boltzmann equation with energy dependent inelasticity. Let us 
emphasize that  it also gives an alternative proof of existence of
solution for the elastic Boltzmann equation with initial datum having 
only finite mass and kinetic energy (but possibly infinite entropy). 

\smallskip
(ii) A generalization of the DiBlasio uniqueness Theorem for the
elastic hard spheres Boltzmann equation (see \cite{DiB,G86,W94,MW99}) and for the
inelastic hard spheres Boltzmann equation with constant normal restitution
coefficients (see \cite{GPV**,FM**}), to the inelastic hard spheres Boltzmann equation 
with energy dependent normal restitution coefficients 
(see also \cite{DM***} where similar tools are developed).

\smallskip\noindent
For points~(i) and~(ii), one of the main ideas of the proof is an
appropriate use of the change of variables $v_* \to v'$ (for fixed
$(v,\sigma)$) and $v \to v'$ (for fixed $(v_*,\sigma)$) in the spirit of the proof of the so-called ``cancelation lemma'' introduced in~\cite{Vill:99} (see also~\cite{ADVW00}).

\smallskip
(iii) An improvement of the result of propagation of exponential
moments for the elastic Boltzmann equation \cite{B97} and for the
inelastic Boltzmann equation \cite{BGP**}, into a result of instantaneous
appearance of exponential moments. This is obtained by combining
estimates from~\cite{BGP**} together with a simple o.d.e. argument 
introduced in the context of the Boltzmann equation in \cite{W97}.

\medskip
Before we explain our results and methods in details, let us introduce
the problem. 

\subsection{A general framework for the collision operator}

We denote by $B$ the rate of occurance of collisions of two particles with
{\it pre-collisional velocities} $\{v,v_*\}$ which gives rise to
post-collisional velocities
$\{v',v'_*\}$. The collision may be schematically written
      \beqn  \label{InelastiColl}
      \{v\} + \{v_*\} \stackrel{B}{\longrightarrow} \{v'\} + \{v'_*\}
      \quad\hbox{with}\quad
       \left\{
       \begin{array}{lcl}
       && \!\!\!\!\!\!\!\!\!\!\!\! v'+v'_* = v+v_*\\ \\
       && \!\!\!\!\!\!\!\!\!\!\!\! |v'|^2+|v'_*|^2 \le |v|^2+|v_*|^2.\\
       \end{array} \right.
      \eeqn
More precisely, for any fixed {\it pre-collisional velocities} $v, \,
v_* \in \R^N$, we introduce a parametrization by $z \in D := \{ w \in
\RR^N; \, |w| \le 1 \}$ of all possible resulting {\it post-collisional
velocities} $\{v',v'_*\}$ in (\ref{InelastiColl}) in the following way:
      \beqn
      \label{Param2}
      \left\{
       \begin{array}{lcl}
       && \!\!\!\!\!\!\!\!\!\!\!\! v' = (v+ v_*)/2 +  z \, |v_*-v| /2
\\
\\
       && \!\!\!\!\!\!\!\!\!\!\!\! v_*' =(v+ v_*)/2 -  z \, |v_*-v| /2.\\
       \end{array}
      \right.
      \eeqn
The {\it collision rate} $B$  takes the form
  \beqn \label{shapeB1}
        B =  |u| \, b, \quad b = \alpha(\EE) \, \beta(\EE,u;dz)
        \eeqn
where  $u=v-v_*$ is the {\it relative velocity}, $\alpha$ is an intensity 
coefficient, $\beta$ is the normalized cross-section (it is a 
probability measure on $D$ for any fixed $\EE,u$), and $\EE$ is the {\it kinetic energy} of
the distribution $f$, defined by
\[
\EE :=  \EE(f) = \int_{\RR^N} f \, |v|^2 \, dv.
\]
The expression (\ref{shapeB1}) reflects the fact that we are dealing
with {\it hard spheres} which undergo contact interactions.  The term
$|u| \, \alpha(\EE)$ corresponds to the rate of collisions of two
particles with  {\it pre-collisional velocities} $v, \, v_* \in \R^N$,
while the term $\beta$ corresponds to the conditional distributional
probability  to obtain the two  {\it post-collisional velocities}
$\{v',v'_*\}$.
The non-negative real $|z|$ is the {\it restitution coefficient} which
measures the loss of energy in the collision, since
\beqn\label{Deltav} |v'|^2 + |v'_*|^2 - |v|^2 - |v_*|^2 = - {1 \over 2}
   (1 - |z|^2)
\, |v_*-v|^2 \le 0.
\eeqn
In the above formula, $|z| = 1$ corresponds to an elastic collision
while $z=0$ corresponds to a completely inelastic collision (or {\it
sticky collision}).

\medskip
The bilinear collision operator $Q(f,f)$ models the interaction of
particles by means of inelastic binary collisions
(preserving mass and total momentum but dissipating kinetic
energy).
More precisely, we define the collision operator by its action on test
functions (which is related to the evolution of the {\it observables} of the probability
density).
Taking $\varphi =\varphi(v)$ to be some well-suited regular
function, we introduce the following weak formulation of the collision
operator (valid under the symmetry assumption~(\ref{shapeB2}) below on $\beta$) 
      \beqn \label{InelasticOp}\qquad
      \langle Q(f,f),\varphi \rangle \, := \,
      {1 \over 2} \int_{\RR^N} \! \int_{\RR^N} \, f_* \, f \int_{D}
      (\varphi'_* + \varphi'- \varphi -
      \varphi_*) \, B(\EE,u;dz) \, dv \, dv_*.
      \eeqn
Here and below we use the shorthand notations $\psi := \psi(v)$, $\psi_*
:= \psi(v_*)$, $\psi' := \psi(v')$ and $\psi'_* :=\psi(v'_*)$
for any function $\psi$ on $\RR^N$.

\medskip
A first simple consequence of the definition of the
operator~(\ref{InelasticOp}) and of the parametrization~(\ref{Param2})
is that mass and momentum are conserved
      \[
      {d \over dt} \int_{\RR^N} f \,
       \left(
        \begin{array}{ll}
        1 \\
        v
        \end{array}
       \right)
      \, dv = 0,
      \]
a fact that we easily derive
(at least formally), multiplying
the equation~(\ref{eqB1}) by $\varphi = 1$ or  $\varphi=v$
and integrating in the velocity variable (using (\ref{InelasticOp})).
In the same way, multiplying
equation~(\ref{eqB1}) by $\varphi = |v|^2$, integrating and using
(\ref{Deltav}) and (\ref{InelasticOp}), we obtain that the 
kinetic energy is dissipated
      \beqn \label{eqdiffEE}
      {d \over dt} \EE(t) =  - D(f) \le 0,
      \eeqn
where we define  the energy dissipation functional $D$ and the energy
dissipation rate $\Delta$, which measures the (averaged) inelasticity of
collisions, by
   \bean\label{defD}
    D (f)&:= & \int_{\RR^N} \! \int_{\RR^N} f \, f_* \, |u|^3 \,
      \Delta(\EE,u) \,\, dv \, dv_*, \\ \label{defbeta}
    \Delta(\EE,u) &:=& \frac14 \, \int_D (1-|z|^2) \, b(\EE,u;dz) \ge 0.
   \eean

\medskip
Finally, we introduce the {\em cooling time}, associated to the process
of cooling (possibly in finite time) of granular gases:
       \beqn  \label{CoolingTime}
       T_c := \inf \Big\{ T \ge 0, \ \EE(t) = 0 \,\, \forall \, t > T \Big\}
         = \sup \Big\{ S \ge 0, \ \EE(t) > 0 \,\, \forall \, t < S \Big\}.
       \eeqn
This cooling effect (or collapse) is one of the main motivations for
the physical and mathematical study of granular media.

\smallskip
The Boltzmann equation~(\ref{eqB1}) is complemented with an initial
condition~(\ref{eqB2})
where the initial datum is supposed to satisfy the moment conditions
      \beqn  \label{initialcond}
      0 \le f_{\mbox{\scriptsize{in}}} \in L^1_q(\RR^N), \qquad
      \int_{\RR^N} f_{\mbox{\scriptsize{in}}} \, dv  = 1, \qquad
      \int_{\RR^N} f_{\mbox{\scriptsize{in}}} \, v \, dv  = 0
      \eeqn
for some $q \ge 2$. Notice that we can assume without loss of generality
the two last moment conditions in~(\ref{initialcond}), since we
may always reduce to that case by a scalling and translation argument.
Here we denote, for any integer $q \in \NN$, the Banach space
      \[
      L^1_q = \left\{f: \RR^N \longrightarrow \RR \hbox{ measurable}; \;
\;
      \| f \|_{L^1_q} := \int_{\RR^N} | f (v) | \, (1+|v|^q) \,  dv
      < \infty \right\}.
      \]
We also define the weighted Sobolev spaces $W^{k,1} _q$
($q \in \R$ and $k \in \N$) by the norm
     \[ \| f \|_{W^{k,1} _q} = \sum_{|s| \le k}
       \|\partial^s f \, (1+|v|^q) \|_{L^1}.  \]
We introduce the space of normalized probability measures on
$\RR^N$, denoted by $M^1 (\RR^N)$, and the space $BV_q(\RR^N)$ ($q \in
\RR$) of (weighted) Bounded Variation functions, defined as the set of the weak limits in
$\DD'(\RR^N)$ of sequences of smooth functions which are bounded in
$W^{1,1}_q(\RR^N)$. Throughout the paper we denote by ``$\mbox{C}$'' 
various constants which do not depend on the collision rate $B$.

\subsection{Mathematical assumptions on the collision rate}

Let us state the basic assumptions on the collision rate $B$:
     \begin{itemize}
        \item The probability measure $\beta$ satisfies the symmetry property
        \beqn  \label{shapeB2}
        \beta(\EE,u;dz) = \beta(\EE,-u;-dz).
        \eeqn
     \item For any $\varphi \in C_c(\R^N)$ the functions
        \beqn \label{shapeB3}
        (v,v_*,\EE) \mapsto \int_D \varphi(v')  \,
        \beta(\EE,u;dz) \quad\hbox{and}\quad \EE \mapsto \alpha(\EE)
        \eeqn
     are continuous on $\R^N \times \R^N \times (0,\infty)$ and
$(0,\infty)$ respectively.
\item The probability measure $\beta$ satisfies the following angular
spreading property: for any $\EE >0$, there is a function $j_\EE(\e) \ge 0$ such that
      \beqn \label{shapeB4}
      \qquad
      \forall \, \e >0, \ u \in \R^N \quad\,\,\,
      \int_{\big\{|\hat{u}\cdot z|\in [-1,1] \setminus [-1+ \e;1-\e]\big\}}
      \beta (\EE,u;dz) \le  j_\EE(\e) 
      \eeqn
and $j_\EE(\e) \to 0$ as $\e \to 0$ uniformly according to $\EE$ when it is restricted to a compact set of
$(0,+\infty)$.
\end{itemize}

We will sometimes restrict our analysis to a kind of {\em generalized (energy dependent) visco-elastic model} 
assuming that the cross-section $b$ reduces to an absolutely continuous measure
according to the Hausdorff measure on  the sphere
\beqn\label{defCue}
     \mathcal{C}_{u,e} = \frac{1-e}{2} \, \hat u + \frac{1+e}{2} \,
\Sph^{N-1}.
\eeqn
More precisely, we assume that 
\beqn
\label{bCue} b(\EE,u;dz)= \delta_{\{z=(1-e)\hat u/2 +
(1+e)\sigma/2 \}} \, \tilde{b}(\EE,|u|,\hat{u}\cdot \sigma) \, d\sigma
\eeqn
where $d \sigma$ is the uniform measure on the unit sphere,
$\tilde b$ is a non-negative measurable function and 
$e : (0,\infty) \times \R^N \times [-1,1] \to [0,1]$, 
$e = e (\EE,|u|,\hat u \cdot \sigma)$ 
is a continuous function. For a vector $x \in \R^N \backslash \{ 0 \}$, we define $\hat{x}=x/|x|$
and $\Sph^{N-1}$ stands for the unit sphere of $\R^N$.
Roughly speaking, the {\em generalized energy dependent visco-elastic model} corresponds 
then to the case where for any direction $\hat z \in \Sph^{N-1}$, 
the post-collisional velocities $(v',v'_*)$ such that $(v' - v'_*)/|v'-v'_*|= \hat z$ 
are uniquely (or deterministically) defined by the pre-collisional velocities $(v,v_*)$.

\smallskip
For the uniqueness of the energy coupled models, we shall need the
following additional assumption:

\begin{itemize}
\item[{\bf H1.}] The cross-section $b$ satisfies (\ref{bCue}) 
with $\tilde b$ bounded, $e=e(\EE)$ and the following
locally Lipschitz conditions holds: for any compact subset $K \subset (0,\infty)$ 
there exists a constant $L_K\in (0,\infty)$ such that  for any $\EE, \EE' \in K$
\beqn\label{H1-1}
  \sup_{u \in \R^N} \|Ê\tilde b (\EE', u, . ) - \tilde b (\EE, u, . ) \|_{L^1(\Sph^{N-1})} 
  \le L_K \, |\EE' - \EE|
\eeqn
and 
\beqn\label{H1-2}
 |Êe (\EE' ) - e (\EE ) |  \le L_K \, |\EE' - \EE|. 
\eeqn
\end{itemize}

\medskip
In the study of the cooling process, we always assume:

\begin{itemize}
\item[{\bf H2.}] The energy dissipation rate $\Delta(\EE,u)$ 
in (\ref{defbeta}) is continuous
on $(0,+\infty) \times \RR^N$ and satisfies
      \beqn \label{hypbeta1}
      \Delta(\EE,u) > 0
      \quad  \forall \, u \in \RR^N, \, \, \EE > 0.
      \eeqn
\end{itemize}

We will also need one of the two following additional assumptions:

\begin{itemize}
\item[{\bf H3.}]  For any $\EE_0,\EE_\infty \in (0,\infty)$ (with
$\EE_0 \ge \EE_\infty$) there exists $\psi$ such that
      \beqn \label{hypbeta2}
      \Delta(\EE,u) \ge \psi(|u|) \quad\forall \, \EE \in
(\EE_\infty,\EE_0),
      \, \,\, \forall \, u \in \R^N,
      \eeqn
with $\psi \in C(\R_+,\R_+)$ and such that for any $R > 0$ there exists
$\psi_R > 0$ with 
\beqn \label{hypbeta3}
\psi(|u|) \ge \psi_R \,  |u|^{-1} \quad \forall \, u \in \R^N, \, \,
|u| > R/2.
\eeqn
This assumption is quite natural. In particular, it holds
for a {\it ``normal'' granular media}.

\item[{\bf H4.}]  The cross-section $b$ satisfies  (\ref{bCue}) 
with $e=e(\EE,|u|)$ and there exists $b_0,b_1 \in (0,\infty)$ such that $b_0 \le \tilde b \le b_1$ a.e. and
$x \mapsto \tilde b(\EE,|u|,x)$ is nondecreasing and convex on $(-1,1)$ for any fixed $\EE \in (0,\infty)$
and $u \in \R^N$.
\end{itemize}

Notice that under assumption (\ref{bCue}) with $\tilde b = \tilde b(\hat{u}\cdot \sigma)$ 
and $e = e(\EE,u)$ the energy dissipation rate just writes
     \beqn\label{betaEu}
     \Delta(\EE,u) = C_N \, (1 - e^2),
     \eeqn
where $C_N$ is a constant depending on the dimension.

\medskip
Let us emphasize that the classical Boltzmann collision operator for
inelastic hard spheres with a constant normal restitution coefficient $e \in
[0,1]$, as studied in~\cite{BCG00} and~\cite{GPV**}, is included as a
particular case of our model, and satisfies all the assumptions above.
But the formalism described from (\ref{InelastiColl}) to (\ref{shapeB4})
is much more general than this case. In particular, we may also
consider:
\smallskip

1. Uniformly inelastic collision processes such that
     \beqn\label{e0}
     \qquad  \exists \, z_0 \in [0,1) \quad\hbox{s.t.}\quad
     \hbox{supp}\, B(\EE,u,.) \subset D(0,z_0) \quad \forall \, u
     \in \RR^N, \,\,\, \forall \, \EE > 0,
     \eeqn
which includes the {\em sticky particles model} when $z_0 = 0$.
\smallskip

2. The physically important case~(\ref{defCue},\ref{bCue}) 
of collisions defined by a normal restitution coefficient $e$ and the cross-section $\tilde{b}$
which possibly depend on $\EE$, $|u|$ and $\hat{u} \cdot \sigma$. In particular
it covers the kind of models studied in \cite{BCG00} (where
$e$ depends on $\EE$, and $\tilde b$ is independent
on $\EE$ and $|u|$). It includes also the important case of 
the {\it visco-elastic hard spheres model}  
where $\tilde b = \tilde b(\hat u \cdot \sigma)$ and 
the normal restitution coefficient depends (smoothly) on the 
normal component of the relative velocity, that is 
$|u| |\hat u - \sigma| /2$ in our notation (see~\cite{BriPoe}).  
\smallskip

3. This formalism also covers multidimensional versions
of the kind of models proposed in \cite{Toscani00}, which
corresponds to the case where $b$ is the product of
a measure depending on $|u|$, $|z|$ and a measure
of $\hat{u} \cdot \hat z$ absolutely continuous
according to the Hausdorff measure. 
One easily checks that our assumptions~(\ref{shapeB1},\ref{shapeB2},\ref{shapeB3},\ref{shapeB4}) 
on the collision rate are quite natural for this kind of models as well.
Note that our measure framework for $B$ can also models situations 
where, in the opposite to the {\em generalized 
visco-elastic case}, there is some stochasticity or uncertainty on the degree of inelasticity 
of the collisions, for instance due to some experimental noise, or due to the fact that particles 
in the gas are a mixture of different inelasticity behaviors, 
which are therefore handled statistically. 
\medskip

   The fact that $b$ is a finite measure on $D$ allows to define
the splitting $Q=Q^+ - Q^-$
where $Q^+$ and $Q^-$ are defined in weak form by
      \beqn \label{defQ+dual}
      \langle Q^+(g,f),\varphi \rangle \, := \, \int_{\RR^N} \!
\int_{\RR^N} \, g_* \, f
      \int_{D} \varphi'  \, |u| \, b(\EE,u;dz) \, dv \, dv_*
      \eeqn
and
      \beqn
      \langle Q^-(g,f),\varphi \rangle \, := \, \int_{\RR^N} \!
\int_{\RR^N} \, g_* \, f
      \int_{D} \varphi \, |u| \, b(\EE,u;dz) \, dv \, dv_*, 
      \eeqn
where $v'$ is defined by~(\ref{Param2}). 
A straightforward computation shows that it is possible to give
a very simple strong form of $Q^-$ as follows
      \beqn
      Q^-(g,f) = L(g) \, f,
      \eeqn
where $L$ is the convolution operator
      \beqn
      L(g)(v) := \alpha(\EE) \, \int_{\RR^N} \, g(v_*) \, |v-v_*| \,
dv_*.
      \eeqn

Under assumption (\ref{bCue}), the expression of $Q^+$ reduces
to
      \beqn
      \langle Q^+(g,f),\varphi \rangle \, := \, \int_{\RR^N} \!
\int_{\RR^N} \, g_* \, f
      \, |u| \int_{\Sph^{N-1}} \varphi'  \, \tilde b(\EE,|u|,
\hat{u}\cdot
\sigma) \, d\sigma \, dv \, dv_*, 
      \eeqn
where $v'$ is defined by the formula (deduced from (\ref{Param2}) and (\ref{defCue}))
\beqn\label{vprimvisco}
 v' = v  -  \frac{1+e}{4} \, \Big[ u - |u| \sigma \Big],
 \quad
 v'_* = v_* +  \frac{1+e}{4} \, \Big[ u - |u| \sigma \Big].
 \eeqn

\subsection{Statement of the main results}

Let us now define the notion of solutions we deal with in this paper.

\begin{defin} \label{def1}
Consider an initial datum $f_{\mbox{\scriptsize{{\em in}}}}$ satisfying
(\ref{initialcond}) with $q=2$.
A nonnegative function $f$ on $[0,T] \times \R^N$ is said to be a
solution to the Boltzmann equation (\ref{eqB1})-(\ref{eqB2}) if
\beqn
\label{fL12}
f \in C([0,T];L^1_2(\R^N)),
\eeqn
and if (\ref{eqB1})-(\ref{eqB2})  holds in the sense of distributions,
that is,
\beqn
\label{dualBeq}
\int_0^T \left\{ \!\int_{\R^N}  f \, {\partial \phi \over \partial t}
dv  \langle Q(f,f),\phi \rangle  \right\} \, dt
= \int_{\R^N} f_{\mbox{\scriptsize{in}}} \, \phi(0,\cdot) \, dv
\eeqn
for any  $\phi \in C^1_c([0,T) \times \R^N)$.
\end{defin}

It is worth mentioning that (\ref{fL12}) ensures that the collision
term $Q(f,f)$ is well defined as a function of $L^1(\R^N)$. Indeed, on
the one hand, we deduce from
$f \in C([0,T];L^1_2(\R^N))$ that $\EE(t) \in K_1$ on $[0,T]$ and thus
$\alpha(\EE(t)) \in K_2 $ on $[0,T]$ for some compact sets $K_i \subset
(0,\infty)$. On the other hand, from the
dual form~(\ref{defQ+dual}) it is immediate that $Q^\pm$ is bounded from
$L^1 _1 \times L^1 _1$ into $L^1$, with bound $\alpha(\EE)$
(see also \cite{GPV**,MM05} for some {\sl strong forms} of the
$Q^+(f,f)$ term). It turns out that a solution $f$, defined as
above, is also a
solution of (\ref{eqB1})-(\ref{eqB2}) in the mild sense:
     \[
     f(t,\cdot) = f_{\mbox{\scriptsize{in}}} + \int_0^t  Q(f(s,\cdot),f(s,\cdot)) \, ds 
    \quad\hbox{a.e. in} \quad \R^N.
     \]
Another straightforward consequence is that if $f \in L^\infty([0,T),L^1_q)$
then $f$ satisfies the following {\it chain rule}
\beqn
\label{renormBeq}
{d \over dt} \int_{\R^N} \Xi(f) \, \phi \, dv = \langle  Q(f,f) ,
\Xi'(f)
\, \phi \, \rangle   \quad\hbox{in} \quad \DD'([0,T)),
\eeqn
for any $\Xi \in C^1(\RR) \cap W^{1,\infty}(\RR)$,  $\phi \in
L^\infty_{1-q}(\R^N)$, in the sense of distribution on $[0,T)$.  

\medskip
Let us state the main results of this paper.
First, we give a Cauchy Theorem valid when the collision rate $B$ is
independent on the kinetic energy.
      \begin{theo}\label{Lcase}
      Assume that $B$ satisfies the assumptions
      (\ref{shapeB1})-(\ref{shapeB2})-(\ref{shapeB3})-(\ref{shapeB4}) 
      with $b = b(u;dz)$: the cross-section does not depend on the kinetic energy.
      Take an initial datum $f_{\mbox{\scriptsize{{\em in}}}}$ satisfying (\ref{initialcond}) with $q= 3$.
      Then
      \begin{itemize}
      \item[(i)]
      For all $T >0$, there exists a unique solution $f \in
      C([0,T];L^1_2) \cap L^\infty (0,T;L^1_3)$ to the Boltzmann equation
      (\ref{eqB1})-(\ref{eqB2}).
      This solution conserves mass and momentum,
        \beqn\label{consmq}
        \int_{\RR^N} f(t,v) \, dv  = 1, \qquad
        \int_{\RR^N} f(t,v) \, v \, dv  = 0 \qquad \forall \, t  \in [0,T],
        \eeqn
      and has a positive and decreasing kinetic energy
        \beqn\label{decE}
        0 < \EE(t_2) \le \EE(t_1) \le \EE_{\mbox{\scriptsize{{\em in}}}} = \EE(0) \qquad
\forall \, t_i \in [0,T],
        \,\,\, t_1 \le t_2.
        \eeqn
      In particular, the life time of the solution (as introduced in~(\ref{CoolingTime})) is
      $T_c = + \infty$.
      \item[(ii)] Moreover,
      assuming {\bf H2-H3} or {\bf H2-H4} (with $e$ and $\tilde b$
      independent on the kinetic energy), there holds
        \beqn  \label{asymptTc}\qquad
        \EE(t) \to 0 \, \hbox{ and } \, f(t,.) \ \wto \  \delta_{v=0}
        \,\hbox{ in }\, M^1(\RR^N)\hbox{-weak}*
        \, \hbox{ when } \, t \to T_c.
        \eeqn
      In other words, the cooling process does not occur in finite time, but
      asymptotically in large time.
      \end{itemize}
      \end{theo}

\begin{rems} \label{remLinear}
Let us discuss the assumptions and
conclusions of this theorem.
\smallskip

1. Under assumption {\bf H4} and when the collision rate is independent
on the
kinetic energy, one can prove in fact that there exists a unique
solution $f \in  C([0,\infty);L^1)$ satisfying (\ref{consmq}) and
(\ref{decE})
for any initial condition $f_{\mbox{\scriptsize{{\em in}}}}$ satisfying
(\ref{initialcond}) with
$q = 2$. The proof is quite
more technical and we
refer to \cite{MW99} where the result is presented for the true elastic
collision Boltzmann equation; nevertheless the proof may be readily
adapted to the inelastic collisional framework.
\smallskip

2. The existence and uniqueness part of Theorem \ref{Lcase} (point (i))
extends to a collision rate $B = B(u;dz) \ge 0 $ which satisfies the
sole assumptions
     \[ \left\{
     \begin{array}{l}
\quad B(-u;-dz) = B(u;dz), \vspace{0.2cm} \\ \displaystyle
\quad \int_D B \, dz \le C_0 \, (1+|v|+|v_*|) \vspace{0.2cm} \\
\displaystyle
\quad (v,v_*)  \mapsto \int_D \varphi(v') \, B(u;dz) \in C(\R^N
\times
     \R^N) \qquad \forall \, \varphi  \in C_c(\R^N)
     \end{array} \right. \]
for some constant $C_0 \in \RR_+$. This corresponds to the so-called
cut-off hard potentials (or variable hard spheres) model 
in the context of inelastic gases.
\smallskip

3. For a uniformly dissipative collision model, {\em i.e.}, such that
      \[
      \Delta (u) \ge \Delta_0 \in (0,\infty),
      \]
a fact which holds under assumption (\ref{e0}) or under assumption {\bf H4}
with a normal restitution coefficient $e$ satisfying $e(|u|) \le e_0 \in
[0,1)$ for any
$u \in \R^N$, we may prove the additionnal {\em a priori} bound
      \[ \int_0^{+\infty} \| f(t,.) \|_{L^1_3} \, dt \le
C\big(\|f_{\mbox{\scriptsize{{\em in}}}}
      \|_{L^1_2},\Delta_0\big). \]
As a consequence, one can easily adapt the proof of existence and
uniqueness in Theorem~\ref{Lcase} and then one can easily
establish that the existence part of Theorem~\ref{Lcase}  holds for
any initial datum $f_{\mbox{\scriptsize{{\em in}}}}$ satisfying
(\ref{initialcond}) with $q = 2$.
\smallskip

4. The existence and uniqueness part of Theorem~\ref{Lcase} (point (i))
immediately extends to a time dependent collision rate
$B = |u| \, \gamma(t) \, b(t,u;dz)$ where $b(t,u;\cdot)$
is a probability measure for any $u \in \RR^N$, $t \in [0,T]$
such that $b(t,u;dz) = b(t,-u;-dz)$, and
$\gamma(t)$ is a non-negative function in $L^\infty(0,T)$.
\smallskip

5. Finally let us emphasize that Theorem~\ref{Lcase} applies to 
the important (non-coupled) model of visco-elastic hard spheres. 
Indeed the collision rate of this model satisfies assumptions 
(\ref{shapeB1},\ref{shapeB2},\ref{shapeB3},\ref{shapeB4}) as 
well as {\bf H2} and {\bf H3}, with $\tilde b$ and $e$ independent 
of $\EE$. We refer to the work in preparation~\cite{MMvisco**} 
which shall be devoted to the detailed study of this particular model. 
\end{rems}

Now, let us turn to the case where the collision rate depends on the
kinetic energy of the solution.

\begin{theo}\label{NLcase}
Assume now that $B$ satisfies the
assumptions~(\ref{shapeB1})-(\ref{shapeB2})-(\ref{shapeB3})-(\ref{shapeB4}) 
and that the cross-section $b = b(\EE,u;dz)$ indeed depends on the
kinetic energy $\EE$. Take an initial datum $f_{\mbox{\scriptsize{{\em
in}}}}$ satisfying
(\ref{initialcond}) with $q = 3$.

\begin{enumerate}

        \item[(i)] There  exists at least one maximal solution $f \in
C([0,T];L^1_2) \cap
        L^\infty(0,T;L^1_3)$, $\forall \, T \in (0,T_c)$, for some $T_c
\in (0,+\infty]$,  which  satisfies
        the conservation laws (\ref{consmq}) and the decay of the
kinetic energy (\ref{decE}).

        \item[(ii)] If the collision rate satisfies the additional assumption {\bf H1}, 
        and the initial datum
        satisfies the additional assumption $f_{\mbox{\scriptsize{{\em
in}}}} \in BV_4 \cap L^1_5$,
        then this solution is unique among the class of functions
$C([0,T],L^1_2) \cap L^\infty (0,T;L^1_3)$, for any $T \in  (0,T_c)$.

        \item[(iii)] The asymptotic convergence (\ref{asymptTc}) holds
under
        the additional assumptions {\bf H2-H3} or {\bf H2-H4}.

        \item[(iv)] If one of following assumptions a. or b. is satsfied, then $T_c=+\infty$: \\ 
        a. $\alpha$ is bounded near $\EE = 0$ and $j_\EE$ converges
        to $0$ as $\e \to 0$ uniformly near $\EE = 0$; \\  
        b. $B$ satifies {\bf H4}, $\Delta$ is bounded by 
        an increasing function $\Delta_0$ which only depends on the energy, and
        $f_{\mbox{\scriptsize{{\em in}}}}\, e^{a_\eta \, |v|^\eta} \in L^1$ with $\eta \in (1,2]$,
        $a_\eta > 0$.

        \item[(v)] If $\Delta(\EE,u) \ge \Delta_0 \, \EE^\delta$ with
$\Delta_0 > 0$
        and $\delta < -1/2$, then $T_c < + \infty$.
\end{enumerate}
\end{theo}

\begin{rem}  Under the assumptions of point (ii) on the initial datum,
by using a bootstrap {\em a posteriori} argument as introduced in
\cite{MW99}, one can  prove that there exists a unique solution
$f \in C([0,\infty);L^1)$ satisfying (\ref{consmq}) and (\ref{decE})
for any initial condition $f_{\mbox{\scriptsize{{\em in}}}}$ satisfying
(\ref{initialcond}) with $q  > 4$ and $f_{\mbox{\scriptsize{\em in}}} \in
BV_4$.
\smallskip
\end{rem}

\subsection{Plan of the paper}

We gather in Section~\ref{sec:Q} some new integrability estimates on the
collision operator which can be of independent interest.
We prove convolution-like estimates in Orlicz spaces for the {\it gain
term}. We give then estimates on the {\it global operator} in Orlicz
space, which show essentially  that even if the bilinear collision
operator is not bounded, its evolution semi-group is bounded in any
Orlicz space (with bound depending on time).
In Section~\ref{sec:cauchy:L} we start looking at {\em solutions} of
the Boltzmann equation. We prove Povzner lemma and several moments
estimates in $L^1$, from which we deduce the existence and uniqueness
part of Theorem~\ref{Lcase}.
In Section~\ref{sec:cauchy:NL}, we extend the existence result to
collision rates depending on the kinetic energy of the solution by
proving a weak stability result on the basis of (local in time)
non-concentration estimates obtained by the study of
Section~\ref{sec:Q}, to
obtain the existence part of Theorem~\ref{NLcase}. The uniqueness part
of Theorem~\ref{NLcase} is obtained by proving a strong stability
result valid for smooth solution.
In Section~\ref{sec:CP} we study the cooling process and prove the
remaining parts of Theorem~\ref{Lcase} and Theorem~\ref{NLcase}.


\section{Estimates in Orlicz spaces}\label{sec:Q}
\setcounter{equation}{0}
\setcounter{theo}{0}


In this section we gather some new functional estimates on the
collision operator in Orlicz spaces, that will be used
in the sequel to obtain (local in time) non-concentration
estimates. Let us introduce the following
decomposition $b = b^t _\e + b^r _\e$
of the cross-section $b$ for $\e \in (0,1)$:
\beqn\label{defbtbr}
      \left\{
        \begin{array}{ll}
        b^t _\e (\EE,u;dz) = b(\EE,u;dz) \, {\bf 1}_{\{-1+\e \le
        \hat{u}\cdot z \le 1 - \e \}} \vspace{0.2cm} \\
        b^r _\e (\EE,u;dz) = b(\EE,u;dz) - b^t _\e (\EE,u;dz)
        \end{array}
      \right.
\eeqn
where ${\bf 1}_{\{-1+\e \le \hat{u}\cdot z  \le 1 - \e \}}$ denotes
the usual indicator function of the set
$\{-1+\e \le \hat{u}\cdot z \le 1 - \e \}$.
When no confusion is possible the subscript $\e$ shall be omitted.

In the sequel, $\Lambda$ denotes a function $C^2$ strictly increasing,
convex
satisfying the assumptions~(\ref{eq:hypLamb1}), (\ref{eq:hypLamb2}) and
(\ref{eq:hypLamb3}) (see the apppendix). This function defines
the Orlicz space $L^\Lambda (\RR^N)$, which is a Banach space
(see the definition in the appendix).

\subsection{Convolution-like estimates on the gain term}

In this subsection we shall prove convolution-like estimates in Orlicz
spaces.
These estimates extend existing results in Lebesgue spaces:
see~\cite{G86,G88,MV**,DM**} in the elastic case and~\cite{GPV**} in
the inelastic case.
The proof relies only upon elementary tools, essentially  Young's
inequality, in the spirit of~\cite{DM**}. 
Moreover it has several advantages: 
its simplicity, the fact that it handles only the dual form
of $Q^+$ and
the fact that it is naturally well-suited to deal with Orlicz spaces,
since it is based on Young's inequality.

As shown by the formula for the differential of the Orlicz norm in the
appendix, the crucial quantity to estimate is
      \[ \int_{\RR^N} Q^+(f,f) \, \Lambda'
\left(\frac{f}{\|f\|_{L^\Lambda}}
         \right) \, dv. \]
Most of the difficulty is related to the fact that the bilinear
operator $Q^+$ is not
bounded because of the term $|v-v_*|$ in the collision rate.
Nevertheless it is
possible to prove a compactness-like estimate with respect to this
algebraic weight.
When combined with the damping effect of the loss term this estimate
shall show that the evolution semi-group of the global collision
operator
is bounded in any Orlicz space.

Let us state the result
\begin{theo}\label{theo:convQ+}
Assume that $B$ satisfies (\ref{shapeB1})-(\ref{shapeB2})-(\ref{shapeB3})-(\ref{shapeB4}). 
For any function $f \in L^1 _1 \cap L^\Lambda$,  for any $\e \in
(0,1)$, there is an
explicit constant $C^+ _\EE(\e)$ such that
\bear\label{eq:convol} \nonumber
     \int_{\RR^N} Q^+(f,f) \, \Lambda'
\left(\frac{f}{\|f\|_{L^\Lambda}} \right) \, dv
     \le \alpha(\EE) \,  \left[  C^+_\EE(\e) \, N ^{\Lambda^*}
     \left( \Lambda' \left( \frac{|f|}{\|f\|_{L^\Lambda}}  \right)
\right)
      \|f\|_{L^1 _1} \|f\|_{L^\Lambda} \right. \\
      +   \left. (2+2^{N+2}) \, j_\EE(\e) \, \|f\|_{L^1 _1} \,
\int_{\RR^N} f \, \Lambda'
       \left(\frac{f}{\|f\|_{L^\Lambda}} \right) \, |v| \, dv \right].
\eear
\end{theo}

\begin{rem} Let us comment on the conclusions of this theorem.

\smallskip
1. We establish estimates for the quadratic
Boltzmann collision operator but similar bilinear estimates could be
proved under additional assumption on $b$, namely that either
{\em no frontal collision occurs}, {\em i.e.}, $b(\EE,u;dz)$
should vanish for $\hat{u}$ close to $z$,
or {\em no grazing collision occurs}, {\em i.e.}, $b(\EE,;dz)$
should vanish for $\hat{u}$ close to $-z$.
For more details on these bilinear estimates and the corresponding
assumptions, we refer to~\cite{MV**} where they are proved
in Lebesgue spaces in the elastic framework.

\smallskip
2. Let us emphasize that for $z \sim 0$ (close to sticky collisions),
the jacobian of the pre-postcollisional change of variable 
$(v,v_*) \to (v',v'_*)$ (both velocities at the same time) is blowing up. However in our
method, we only use the changes of variable $v\to v'$ and
$v_* \to v'$, keeping the other velocity unchanged, and the
jacobians of these changes of variable
remain uniformly bounded as $z \to 0$. This explains
why our bounds includes the sticky particules model,
and are uniform as $z \to 0$.

\smallskip
3. When $\Lambda(t) = t^p/p$, estimate (\ref{eq:convolLp}) just writes
\bear\label{eq:convolLp}
   \qquad  \int_{\RR^N} Q^+(f,f) \, f^{p-1} \, dv \le  \tilde
C^+_\EE(\e) \,  \|f\|_{L^1_1} \, \|f\|_{L^p}^p
     +   \tilde j_\EE(\e) \, \|f\|_{L^1 _1} \,  \|f \, |v|^{1/p}
\|_{L^p}^p,
\eear
for any $\eps \in (0,1)$ and for some explicit constants $\tilde
C^+_\EE(\e)$, $\tilde j_\EE(\e) \in (0,\infty)$
with $ \tilde j_\EE(\e) \to 0$ when $\e \to 0$. Although the quantities
involved in these previously mentioned papers are slightly different, one can see
that estimate (\ref{eq:convolLp}) (or the $L^p$ version of 
Theorem~\ref{theo:OrlQ}) generalizes \cite[Proposition
2.5]{DM**} to the inelastic Boltzmann operator and that  it improves
\cite[Lemma 4.1]{GPV**} because of the better control of the norm $ \|f
\, |v|^{1/p} \|_{L^p}$.

\end{rem}

Let us start with an elementary geometrical lemma that 
we shall need several times in the sequel, 
in order to justify the change of variables $v_* \to v'$ (keeping $v,z$ fixed) 
and $v \to v'$ (keeping $v_*,z$ fixed). 
This lemma is close to the spirit of the proof of these changes of variables 
in the proof of the so-called ``cancellation lemma'' in \cite{Vill:99,ADVW00}. 

\begin{lem}\label{utow} 
For any $z \in D$ and $\gamma \in (-1,1)$ we define the map
\beqn\label{utow1}
\Phi_{z} : \R^N \to \R^N, \quad u \mapsto w = \Phi_{z} (u) := u + |u| \, z,
\eeqn
its Jacobian function $J_{z} := \mbox{{\em det}} \, (D \,  \Phi_z)$ 
and the cone 
$\Omega_\gamma := \{ u \in \R^N\backslash\{0\}, \,\, \hat u \cdot \hat z > \gamma\}$.
Then $\Phi_z$ is a $C^\infty$-diffeomorphism from $\Omega_\gamma$ onto 
$\Omega_{\delta}$ with 
  \[ \delta = \frac{\gamma+|z|}{(1+2\gamma |z| + |z|^2)^{1/2}} \] 
and there exists $C_\gamma \in (0,\infty)$ such that 
\beqn\label{utow2}
C_\gamma^{-1} \le J_z \le C_\gamma \quad\hbox{on}\quad \Omega_\gamma
\eeqn
uniformly with respect to the parameter $z \in D$.
\end{lem}

\smallskip\noindent
{\sl Proof of Lemma \ref{utow}.} 
We may assume $z \not = 0$ since otherwise the conclusion is clear. 
Let start proving that $\Phi_z$ is one-to-one  
on $\Omega_{-1} = \R^N \backslash (\R_- z)$. 
For any $x \in \R^N$ we introduce the decomposition 
$x = x_1 \, \hat z + x_2 :=(x_1,x_2)$ such that $x_1 \in \R$, $x_2 \in \R^N$, 
$x_2 \cdot \hat z = 0$. The expression (\ref{utow1}) then writes equivalently
$$
w_1 = u_1 + (u_1^2 + |u_2|^2)^{1/2} \, |z|, \quad w_2 = u_2.
$$
For any $u,u' \in \Omega_{-1}$ the relation $\Phi_z(u) = \Phi_z(u') =: w$ 
implies immediately $u_2 = u'_2 = w_2$ and we conclude observing 
that for any $z \in D$ and $w_2 \in \R^N$ the map
$$
\varphi_{w_2,|z|} : u_1 \mapsto w_1 := u_1 + (u_1^2 + |w_2|^2)^{1/2} \, |z|
$$
is strictly increasing from $\R$ onto $\R$ if $|z| < 1$, from $\R$ 
onto $\R_+$ if $|z| = 1$ and $w_2 \not=0$, and from $\R_+$ onto $\R_+$ if $|z| = 1$ 
and $w_2=0$. That proves that $\Phi$ is one-to-one. 
Moreover, any point $\hat u = (u_1, u_2) \in \Sph^{N-1}$ 
such that $\hat u_1 = \gamma$ is mapped to the point $w = (\gamma+|z|, u_2)$ 
with square norm  $|w|^2 = 1 + 2 \, \gamma \, |z| + |z|^2$. 
We conclude that $\Phi(\Omega_\gamma) = \Omega_\delta$ thanks to 
the homogeneity property $\Phi_z(r \, u) = r \, \Phi_z(u)$ for any $r > 0$ and $u \in \R^N$. 
We next compute $D\Phi_z(u) =  Id + \hat u \otimes z$ and 
thus $J_z (u) = 1+  \hat u \cdot z$ from which (\ref{utow2}) easily follows. 
Finally, the fact that $\Phi_z$ is a $C^\infty$-diffeomorphism is a direct 
consequence of the local inversion Theorem.  \qed

\medskip\noindent
{\sl Proof of Theorem~\ref{theo:convQ+}.}
Let us denote
\[ \varphi (f) = \Lambda' \left(\frac{f}{\|f\|_{L^\Lambda}} \right).\]
Using the decomposition $b =b^t + b^r$, we control separately the two
terms $I^t$ and $I^r$ in the decomposition
     \bean
     \int_{\RR^N} Q^+(f,f) \, \varphi(f) \, dv
          &= & \int_{\RR^N \times \RR^N \times D} f f_* \varphi(f') \,
|u| \, b^t (\EE,u;dz) \, dv \,
     dv_* \\
           &&+  \int_{\RR^N \times \RR^N \times D} f f_* \varphi(f') \,
|u| \, b^r (\EE,u;dz) \, dv \,
     dv_* =: I^t + I^r.
     \eean
Using the bound
     \[ |u| = |v-v_*| \le |v| + |v_*| \]
we have
     \bean I^t & \le& \int_{\RR^N \times \RR^N \times D} (f|v|) f_*
\varphi(f') \, b^t (\EE,u;dz) \, dv
     \, dv_* \\
                   && + \int_{\RR^N \times \RR^N \times D} f (f_*|v_*|)
\varphi(f')\, |u| \, b^t
     (\EE,u;dz) \, dv \, dv_* =: I^t_1 + I^t_2.
     \eean
For the term $I^t_1$, by applying the Young's inequality (\ref{YoungIneg})
      \[ f_* \varphi(f')  = \|f\|_{L^\Lambda} \left(
\frac{f_*}{\|f\|_{L^\Lambda}} \right) \, \varphi(f')
                            \le \|f\|_{L^\Lambda}  \,
\Lambda\left(\frac{f_*}{\|f\|_{L^\Lambda}}\right)
                                 + \|f\|_{L^\Lambda} \, \Lambda^*
(\varphi(f')),\]
we get 
      \bean
      I^t_1 &\le& \|f\|_{L^\Lambda}  \, \int_{\RR^N \times \RR^N \times
D}
          f |v| \Lambda\left(\frac{f_*}{\|f\|_{L^\Lambda}}\right) \,
             b^t (\EE,u;dz) \, dv \, dv_*  \\
                   && + \|f\|_{L^\Lambda}  \, \int_{\RR^N \times \RR^N
\times D} f |v| \Lambda^*
(\varphi(f'))\, b^t (\EE,u;dz) \, dv \, dv_*
                    =: I^t_{1,1} + I^t_{1,2}.
      \eean
On the one hand, using
      \[ \forall \, x \in \R_+, \quad \Lambda(x) \le x \, \Lambda'(x), \]
which is a trivial consequence of the fact that $\Lambda(0)=0$ and
$\Lambda'$ is increasing, we have
      \[ I^t_{1,1} \le \alpha(\EE) \,  \|f\|_{L^1 _1} \, \int_{\RR^N} f
\,
\varphi(f) \, dv. \]
H\"{o}lder's inequality in Orlicz spaces (\ref{eq:holdOrl}) recalled
in the appendix then yields
      \beqn \label{I111}
       I^t_{1,1}   \le \alpha(\EE) \,
       N ^{\Lambda^*} \left( \Lambda' \left(
\frac{|f|}{\|f\|_{L^\Lambda}}
\right) \right)
        \|f\|_{L^1 _1} \|f\|_{L^\Lambda}.
\eeqn

\smallskip\noindent
On the other hand, using that $\Lambda^*(y) = y \, (\Lambda')^{-1}(y) -
\Lambda( (\Lambda')^{-1}(y))$, we get
      \[  I^t_{1,2}
         \le \int_{\RR^N \times \RR^N \times D} f |v|  \, \varphi(f') f'
\, b^t (\EE,u;dz) \, dv \, dv_*.
    \]
We make the change of variables $v_* \to v'$ (while the other integration 
variables are kept fixed) or more precisely 
$\Psi : (v,v_*,z) \to (v,\psi_{v,z}(v_*),z)$ with 
$\psi_{v,z}(v_*) = v' = v + 2^{-1} \, \Phi_z(v_*-v)$.  
Thanks to the truncation (\ref{defbtbr}) on $b^t_\eps$ and Lemma~\ref{utow}, 
the application $\Psi$ is a $C^\infty$-diffeomorphism from 
$\{ (v,v_*,z) \in \R^{2N} \times D, \,\, \hat u \cdot z \not= 1\}$ onto its image 
and its jacobian $J_\Psi = 2^{-N} \, ( 1 - \hat u \cdot z)$ satisfies 
$| J^{-1}_\Psi |Ê\le 2^N\, \eps^{-1}$ 
on $\{  (v,v_*,z) \in \R^{2N} \times D, \,\, \hat u \cdot z \le 1 - \eps \}$. 
We then get
\bean
I^t_{1,2} 
&\le& \int_{\RR^N \times \R^N \times D} f |v| f' \, \varphi(f') \, 
J^{-1}_\Psi \, b^t (\EE,v-\psi_{v,z}^{-1}(v');dz) \, dv \, dv'  \\
&\le& \alpha(\EE) \, 2^N \e^{-1} \, \|f\|_{L^1 _1}  \, \int_{\RR^N} f \, \varphi(f) \, dv.
\eean
As previously, H\"{o}lder's inequality (\ref{eq:holdOrl}) then yields
      \beqn \label{I112}
       I^t_{1,2}
         \le \alpha(\EE) \, 2^N \e^{-1} \, \|f\|_{L^1 _1} \,
         N ^{\Lambda^*} \left( \Lambda' \left(
\frac{|f|}{\|f\|_{L^\Lambda}} \right) \right)
         \,  \|f\|_{L^\Lambda}.
      \eeqn

Next, the term $I^t_2$ is exactly similar to $I^t_1$, except that one
has to use the change of variable $v \to v' = v_* + 2^{-1} \, \Phi_z(v-v_*)$ instead of $v_* \to v'$.
Therefore,  gathering (\ref{I111}), (\ref{I112}) and
the same estimate for   $I^t_2$, we obtain
     \beqn \label{I1}
     I^t \le 2 \, \alpha(\EE) \, (1+ 2^N \e^{-1}) \, \|f\|_{L^1 _1} \,
           \left[ N ^{\Lambda^*} \left( \Lambda' \left(
     \frac{|f|}{\|f\|_{L^\Lambda}} \right) \right) \right]
         \,  \|f\|_{L^\Lambda}.
     \eeqn

\smallskip\noindent
Finally, for the term $I^r$, we can split it as
\bean
I^r  
&\le&  \int_{\RR^N \times \RR^N \times D} f \, f_* \, \varphi(f')  
            \, {\bf 1}_{\{\hat{u} \cdot z \le 0\}} \, |u| \, b^r (\EE,u;dz) \, dv \, dv_* \\
&& +  \int_{\RR^N \times \RR^N \times D} f \, f_*  \, \varphi(f')  \, {\bf 1}_{\{\hat{u}
     \cdot z \ge 0\}} \, |u| \, b^r (\EE,u;dz) \,  dv \, dv_* =:  I^r_1 + I^r_2.
\eean
For $I^r _1$, we use Young's inequality (\ref{YoungIneg}) on 
$x = f_*$ and $y=\varphi(f')$ to obtain
\bean 
I^r _1 \le \int_{\RR^N \times \RR^N \times D} f \, f_* \, \varphi(f_*)  \, {\bf 1}_{\{\hat{u} \cdot
     z \le 0\}} \,  |u| \, b^r (\EE,u;dz) \, dv \, dv_* \\
          + \int_{\RR^N \times \RR^N \times D} f f' \, \varphi(f')  \, {\bf 1}_{\{\hat{u} \cdot
     z \le 0\}} \, |u| \,  b^r (\EE,u;dz) \, dv \, dv_*.
\eean
In the second integral we make again the change of variable defined by $\Psi$ for which there holds 
$| J^{-1}_\Psi |Ê\le 2^N$ on the domain of integration because of the 
truncation $\hat{u} \cdot z \le 0$. We also observe thanks to a direct 
computation starting from (\ref{InelastiColl}) that under the 
truncation $\hat{u} \cdot z \le 0$ there holds
\[ 
|v-v_*| \le 2 |v'-v| \le 2 (1+|v'|) (1+|v|). 
\]
Hence we obtain
\bean 
I^r _1 
&\le&  (1+2^{N+1}) \, \left( \sup_{u \in \R^N} \int_D b^r(\EE,u;dz) \right) \|f\|_{L^1 _1} \,  \, \int_{\RR^N} f \,
        \varphi(f) \, (1+|v|) \, dv \\ &\le& (1+2^{N+1}) \, \alpha(\EE) \, j_\EE(\e) \, \|f\|_{L^1 _1} \,  \, \int_{\RR^N} 
       f \, \varphi(f) \, (1+|v|) \, dv.
\eean
The term $I^r _2$ is treated similarly using Young's inequality (this time on $x=f$ and $y=\varphi(f')$)  and the change of variable  $v \to v'$ instead of $v_* \to v'$.  It satisfies therefore
the same estimate. Thus we obtain the estimate
\beqn\label{I2}
I^r \le (2+2^{N+2}) \, \alpha(\EE) \, j_\EE (\e) \, \|f\|_{L^1 _1} \,
     \int_{\RR^N} f \, \varphi(f) \, (1+|v|) \, dv.
\eeqn
Defining
\beqn\label{defC+}
C ^+ _\EE (\e) = 2 ( 1 + 2^N \e^{-1} ) + (2+2^{N+2}) \, j_\EE
(\e),\eeqn
we conclude the proof gathering (\ref{I1}) and (\ref{I2}). \qed

\subsection{Minoration of the loss term}

In this subsection we recall a well-known result about
the minoration of the loss term $Q^-$. Let us recall first the following
classical estimate.

\begin{lem}\label{minoLk}
      For any non-negative measurable function $f$ such that
\beqn\label{hyplemQ-}
       f \in L^1_1(\R^N), \qquad \int_{\R^N} f \, dv =1, \qquad
\int_{\R^N} f  \, v\,  dv =0,
\eeqn
we have
\[
    \forall \, v \in \RR^N, \ \ \  \int_{\RR^N} f_* \, |v-v_*| \, dv_*
\ge |v|.
\]
\end{lem}

\medskip\noindent
{\sl Proof of Lemma \ref{minoLk}}. Use Jensen's inequality
      \[ \int_{\RR^N} \varphi(g_*) \, d\mu_* \ge \varphi
\left(\int_{\RR^N}
g_* \, d\mu_* \right) \]
with the probability measure $d\mu_* = f_* \, dv_*$, the measurable
function $v_* \mapsto g_* = v - v_*$
and the convex function $\varphi (s) = |s|$. \qed

\medskip
Then the  proof of the following proposition is straightforward:
      \begin{prop}\label{prop:minQ-}
      Assume that $B$ satisfies~(\ref{shapeB1}). 
      For a non-negative function  $f$ satisfying (\ref{hyplemQ-}), we have
        \begin{equation}\label{eq:minOrl}
        \int_{\RR^N} Q^- (f,f) \, \Lambda'
\left(\frac{f}{\|f\|_{L^\Lambda}}
        \right) \, dv
        \ge \alpha(\EE) \, \int_{\RR^N} f \, \Lambda'
        \left(\frac{f}{\|f\|_{L^\Lambda}} \right) \, |v| \, dv.
        \end{equation}
      \end{prop}


\subsection{Estimate on the global collision operator and {\it a priori}
estimate on the solutions}

Combining Theorem~\ref{theo:convQ+} and Proposition~\ref{prop:minQ-} we
get
\begin{theo} \label{theo:OrlQ}
     Assume that $B$ satisfies (\ref{shapeB1})-(\ref{shapeB2})-(\ref{shapeB3})-(\ref{shapeB4}). 
    Let us consider a non-negative function $f$ satisfying
     (\ref{hyplemQ-}). Then there is an explicit constant $C_\EE$
depending on the
     collision rate
     through the functions $\alpha$ and $j_\EE$ such that
        \[
        \int_{\RR^N} Q(f,f) \, \Lambda' \left(\frac{f}{\|f\|_{L^\Lambda}}
\right) \, dv
        \le C_\EE \,  \left[ N ^{\Lambda^*} \left( \Lambda' \left(
        \frac{|f|}{\|f\|_{L^\Lambda}} \right) \right) \right]
        \|f\|_{L^1 _1} \|f\|_{L^\Lambda}.
        \]
      More precisely, $ C_\EE = \alpha(\EE) \, C^+_\EE (\e_0)$, with
      $\e_0 \mbox{ such that } j_\EE(\e_0) \le (2+2^{N+2})^{-1} \,
\|f\|_{L^1 _1} ^{-1}$ and where $C^+_\EE $ is defined in (\ref{defC+}).
\end{theo}

\medskip\noindent
{\sl Proof of Theorem~\ref{theo:OrlQ}}.
One just has to
combine~(\ref{eq:convol}) and~(\ref{eq:minOrl})
and pick a $\e_0$ small enough such that
      \[ (2+2^{N+2}) \, \|f\|_{L^1 _1} \, j_\EE(\e_0) \le 1. \]
\qed

\begin{cor} \label{cor:Orl}
Assume that $B$ satisfies (\ref{shapeB1})-(\ref{shapeB2})-(\ref{shapeB3})-(\ref{shapeB4}) 
and let us consider a solution $f \in
C([0,T];L^1_2)$ to the Boltzmann equation (\ref{eqB1})-(\ref{eqB2})
associated to an initial datum $f_{\mbox{\scriptsize{{\em in}}}} \in L^1_2$
and to the collision
rate $B$. Assume moreover that  (\ref{consmq}) holds and there exists a
compact set $K \subset (0,+\infty)$  such that
\[
\forall \, t \in [0,T], \quad \EE(t) \in K.
\]
Then, there exists a  $C^2$, strictly increasing and convex function
$\Lambda$ satisfying the   assumptions~(\ref{eq:hypLamb1}),
(\ref{eq:hypLamb2}) and (\ref{eq:hypLamb3}) (which only depends on
$f_{\mbox{\scriptsize{in}}}$) and a constant $C_T$ (which depends on
$K$, $T$ and $B$) such
that
\[ 
\sup_{[0,T]} \|f(t, .) \|_{L^\Lambda} \le C_T.
\]
\end{cor}

\begin{rem}
Let us emphasize that these non-concentration bounds are valid for the
sticky particules
model (in this case they provide an exponentially growing bound in
$L^\Lambda$ for all
times). As a particular case we deduce some explicit bounds on the
entropy when
it is finite initially.
Moreover, since our bounds are uniform as $b \rightharpoonup
\delta_{z=0}$,
we also deduce a proof of the sticky particules limit
(for a cross-section being a diffuse measure
converging to a Dirac mass at $z=0$) by the Dunford-Pettis Lemma.
This shows moreover that this limit is not singular.
\end{rem}

\medskip\noindent{\sl Proof of Corollary~\ref{cor:Orl}. }
Since $f_{\mbox{\scriptsize{in}}} \in L^1(\R^N)$, as recalled in the
appendix,
a refined version of the De la Vall\'ee-Poussin Theorem
\cite[Proposition~I.1.1]{Le}
(see also \cite{LMb,LM02}) guarantees that
there exists a function $\Lambda$ satisfying the properties listed in
the statement of Corollary~\ref{cor:Orl} and such that
     \[ \int_{\RR^N} \Lambda(|f_{\mbox{\scriptsize{in}}}|) \, dv <
+\infty.\]
Then the $L^\Lambda$ norm of $f$ satisfies
     \[ \frac{d}{dt} \|f_t\|_{L^\Lambda}  = \left[ N ^{\Lambda^*} \left(
\Lambda' \left(
     \frac{|f|}{\|f\|_{L^\Lambda}} \right) \right) \right]^{-1}
                        \int_{\RR^N} Q(f,f) \, \Lambda'\left(
     \frac{|f|}{\|f\|_{L^\Lambda}} \right) \, dv
     \]
thanks to Theorem \ref{theo:diffOrlNorm}, and thus using
Theorem~\ref{theo:OrlQ}, we get
     \[ \forall \, t \in [0,T], \quad
     \frac{d}{dt} \|f_t\|_{L^\Lambda} \le C_{\EE(t)} \, \|f_t\|_{L^1 _1}
\,
     \|f_t\|_{L^\Lambda}. \]
Thanks to the assumptions on $B$, the
constant $C_{\EE(t)}$ provided by Theorem~\ref{theo:OrlQ} is uniform
when the
kinetic energy belongs to a compact set. Thus we deduce
    \beqn\label{dnormLLambdadt}
     \forall \, t \in [0,T], \quad
     \frac{d}{dt} \|f_t\|_{L^\Lambda} \le C_K \, \|f_t\|_{L^1 _1} \,
     \|f_t\|_{L^\Lambda}. \eeqn
for some explicit constant $C_K >0$ depending on $K$ and the collision
rate.
We conclude thanks to a Gronwall argument. \qed


\section{Proof of the Cauchy theorem for non-coupled collision
rate}\label{sec:cauchy:L}
\setcounter{equation}{0}
\setcounter{theo}{0}


In this section we fix $T_* > 0$ and we assume that the collision rate
$B$ satisfies
      \beqn \label{shapeB1t}
      B = B(t,u;dz) = |u| \, \gamma(t) \, b(t,u;dz),
      \eeqn
where $b$ is a probability measure on $D$ for any $t \in [0,T_*]$ and
$u \in \RR^N$ satisfying
      \beqn \label{shapeB2t}
      \forall \, t \in [0,T_*], \ \forall \, u \in \R^N, \ \ \
      b(t,u;dz) = b(t,-u;-dz)
      \eeqn
and where $\gamma$ satisfies
      \beqn \label{shapeB3t}
0 \le  \gamma(t) \le \gamma_* \quad \hbox{on} \quad (0,T_*).
      \eeqn

\subsection{Propagation of moments}

In this subsection we establish several moments estimates which are
well known for
the Boltzmann equation with elastic collision, see \cite{B97,MW99,Lu99}
and the references therein,
as well as the recent works \cite{GPV**,BGP**} for the inelastic case.
Let us emphasize that these moment estimates are uniform with respect to
the normal restitution coefficient $e$ or more generally to the
support of $b(t,u;\cdot)$ in $D$.

\smallskip
First we give a result of propagation of  moments valid
for  general collision rates using a
rough version of the Povzner inequality.

\begin{prop}\label{propboundL3}
Assume that $B$ satisfies~(\ref{shapeB1t})--(\ref{shapeB3t}). For any
$0 \le f_{\mbox{\scriptsize{{\em in}}}} \in L^1_q(\R^N)$ with $q > 2$
and $T>0$, there exists $C_T$ such
that any solution $f$ to the inelastic Boltzmann
equation~(\ref{eqB1})-(\ref{eqB2})
on $[0,T]$ satisfies, at least formally,
        \[
  \sup_{[0,T]} \| f(t,\cdot) \|_{L^1 _q} \le C_T.
        \]
\end{prop}

\smallskip\noindent
{\sl Proof of Proposition \ref{propboundL3}.}
We write the proof for the third moment, the general moment estimate
being similar. For any function
$\Psi : \RR^N \to \RR_+$ such that $\Psi (v):= \psi(|v|^2)$
for some function $\psi : \RR_+ \to \RR_+$, the evolution of the
associated  moment is given by
      \[ {d \over dt} \int_{\RR^N} f \, \Psi \, dv
         = \int_{\RR^N \times \RR^N} f \, f_* \, K_\Psi \, dv \, dv_*,
\]
where
      \[ K_\Psi := {1 \over 2}\int_D (\Psi' + \Psi'_* - \Psi - \Psi_*) \,
B(t,u;dz). \]
For $\psi(z) = z^s$, $s > 1$, the function $\psi$ is super-additive,
that is $\psi(x) + \psi(y) \le \psi(x+y)$,
and it is an increasing function. As a consequence,
      \bean
      \Psi' + \Psi'_* - \Psi - \Psi_* &\le& \psi(|v'|^2) + \psi(|v'_*|^2)
-
      \psi(|v'|^2+ |v'_*|^2 ) \\
      && + \psi(|v|^2+ |v_*|^2 ) - \psi(|v|^2 ) - \psi(|v_*|^2 ) \\
      &\le& \psi(|v|^2+ |v_*|^2 ) - \psi(|v|^2 ) - \psi(|v_*|^2 ),
      \eean
which implies
      \[ K_\Psi \le {\gamma(t) \over 2} \, |v-v_*| \, \Big[ \psi(|v|^2+ |v_*|^2 )
      - \psi(|v|^2 ) - \psi(|v_*|^2 ) \Big]. \]
Making the choice $\psi(x) = x^{3/2}$ and using the inequality
     \bear\label{povzner1} \nonumber
     \qquad  (x^{1/2} + y^{1/2}) \, [(x+y)^{3/2} - x^{3/2} - y^{3/2}] &\le&
      C \,  (x^{1/2} + y^{1/2})  \, (x^{1/2}y + xy^{1/2}) \\
      &\le& C \, (2 xy+x^{1/2}y^{3/2} + x^{3/2} y^{1/2})
\eear
for any $x,y > 0$, we get
\beqn\label{Y3gal}
     {d \over dt} \int_{\RR^N} f \, |v|^3 \, dv \le C \, \gamma(t) \,
\int_{\RR^N \times \RR^N}
     f \, f_* \, (|v|^2 \, |v_*|^2  + |v|\, |v_*|^3) \, dv \, dv_*,
\eeqn
and we conclude thanks to a Gronwall argument. \qed

\smallskip
Finally we give a much more precise result on the evolution of moments in
the case when assumption {\bf H4} is made. On the one hand, we state
uniform in time propagation of algebraic moments (as introduced in
\cite{Po62,Ar72,El83}) and exponential moments (for which the first results 
were obtained in \cite{B97}). On the other hand, we prove appearance of
some exponential moments (while appearance of algebraic moments 
was initiated in \cite{Des93,W94,W97})
using carefully estimates developed in \cite{BGP**}.
These estimates may be seen as {\it a priori} bounds, but in fact, by
the
bootstrap argument introduced in \cite{MW99},
they can be obtained {\it a posteriori} for any  solution given by the
existence part of Theorem \ref{Lcase} and Theorem \ref{NLcase}.

\begin{prop}\label{Y3uniform}
We make the assumption {\bf H4} on $B$. A solution $f$ to the inelastic
Boltzman equation~(\ref{eqB1})-(\ref{eqB2}) on $[0,T_c)$ satisfies the
additional
moment properties:
    \begin{enumerate}
     \item[(i)]
For any $s > 2$, there exists $C_s > 0$  such that
\beqn\label{MomentUnif1}
\sup_{t \in [0,T_c)} \|f(t,.) \|_{L^1_s} \le \max\big\{
\|f_{\mbox{\scriptsize{{\em in}}}} \|_{L^1_s},C_s \big\}.
\eeqn
     \item[(ii)]
If $f_{\mbox{\scriptsize{{\em in}}}} \, e^{r \, |v|^\eta} \in
L^1(\R^N)$ for $r>0$ and $\eta \in (0,2]$,
there exists $C_1, r' > 0$, such that
\beqn\label{MomentUnif2}
\sup_{t \in [0,T_c)} \int_{\R^N} f(t,v) \, e^{r' \, |v|^\eta}  \, dv
\le C_1.
\eeqn
       \item[(iii)]
For any  $\eta \in (0,1/2)$ and $\tau \in (0,T_c)$ there exists
$a_\eta,C_\eta \in (0,\infty)$
such that
\beqn\label{expeta}
\sup_{t \in [\tau,T_c)} \int_{\R^N} f(t,v) \, e^{a_\eta \, |v|^\eta} \,
dv \le C_\eta.
\eeqn
\end{enumerate}
Let us emphasize that none of these constants depends on the
inelasticity coefficient $e$ (so that the estimates are uniform with
respect to the  inelasticity of the Boltzmann operator) and that the
constant $C_s, a_\eta, C_\eta$ may depend on
$f_{\mbox{\scriptsize{{\em in}}}}$ only  through
its kinetic energy $\EE_{\mbox{\scriptsize{{\em in}}}}$.
\end{prop}

\begin{rem} The proof of (i) is very classical for the elastic
Boltzmann equation \cite{Po62,Ar72,El83} and it has been extended to
the inelastic operator in \cite{GPV**}. Estimate (ii) has been proved
in \cite{B97} for  the elastic Boltzmann equation and it has been
generalized in \cite{BGP**} to the (stationary) inelastic Boltzmann
equation. We refer  to \cite{B97,MW99,Lu99,CVhand}  for development
around the Povzner inequalities.
Since (ii) is a straightforward consequence of the Povzner inequality
proved in \cite{BGP**}, we just have to prove (iii). Nevertheless,
since the proof of (iii) requires some tools and notations introduced in
\cite{GPV**,BGP**} we begin (step 1 and step 2) by briefly presenting
the proof of (ii). 
Let us emphasize again that (iii) is new even for the elastic equation. 
In the elastic framework, an extension of (iii) to hard potentials 
with cutoff has been used recently in the proof of the
exponential return to equilibrium with explicit rate for initial data
with finite mass and energy, see~\cite{GM:04}.

\end{rem}

\smallskip\noindent
{\sl Proof of Proposition \ref{Y3uniform}. }
Let us define
      \[ m_p := \int_{\R^N} f \, |v|^{2p} \, dv. \]
\medskip\noindent{\sl Step 1. Differential inequalities on the moments.}
Taking $\psi(x) = x^{p/2}$ and $B$ of the above form, there holds
\beqn \label{mp1}
      {d \over dt} m_p = \int_{\R^N} Q(f,f) \, |v|^{2p} \, dv
      = \alpha(\EE) \, \int_{\RR^N \times \RR^N} f \, f_* \, |v-v_*| \,
K_p (v,v_*) \, dv \, dv_*,
\eeqn
where
\beqn \label{mp2}
      K_p(v,v_*) := {1 \over 2}\int_{\Sph^{N-1}} (|v'|^{2p} +
|v'_*|^{2p}  -|v|^{2p} -
      |v_*|^{2p}) \, {\tilde b(\EE,|u|,\sigma\cdot \hat{u}) \over
\alpha(\EE)} \, d\sigma.
\eeqn
From~\cite[Lemma 1, Corollary 3]{BGP**}, there holds
      \beqn \label{mp3}
      K_p(v,v_*) \le \gamma_p \, (|v|^2 + |v_*|^2 )^p - |v|^{2p} -
|v_*|^{2p}
      \eeqn
where $(\gamma_p)_{p=3/2,2,...}$ is a decreasing sequence of real
numbers such that
      \beqn \label{mp4}
      0 <  \gamma_p < \min \left\{1, {4 \over p+1} \right\}
      \eeqn
(notice that the assumptions~\cite[(2.11)-(2.12)-(2.13)]{BGP**}
are satisfied under our assumptions on the collision rate).
Let us emphasize that the estimate~(\ref{mp3}) does not depend on the
inelasticity coefficient $e(\EE,|u|)$.
Then, from~\cite[Lemma~2 and Lemma~3]{BGP**}, we have
      \beqn \label{mp5}
      {1 \over \alpha(\EE)} \int_{\R^N} Q(f,f) \, |v|^{2p} \, dv
      \le \gamma_p \, S_p - (1-\gamma_p) \, m_{p+1/2}
      \eeqn
with
      \[ S_p := \sum_{k=1}^{k_p} \pmatrix{p \cr k} (m_{k+1/2} \, m_{p-k}
+
m_{k}
      \, m_{p-k+1/2}), \]
where $k_p := [(p+1)/2]$ is the integer part of $(p+1)/2$ and
$\pmatrix{p \cr k}$ stands for the binomial coefficient.
Gathering~(\ref{mp1}) and~(\ref{mp5}), we get
      \beqn \label{mp6}\qquad\quad
      {d \over dt} m_p \le \alpha(\EE) \, (\gamma_p \, S_p - (1-\gamma_p)
\,
      m_{p+1/2}) \qquad \forall \, p  = 3/2,2, \dots
      \eeqn
By H\"older's inequality and the conservation of mass,
      \[ m_p^{1+{1 \over 2p}} \le m_{p+1/2} \]
and, by~\cite[Lemma~4]{BGP**}, for any $a \ge 1$, there exists $A > 0$
such that
      \[ S_p \le A \, \Gamma (a\,p+a/2 + 1) \, Z_p \]
with
      \[ Z_p := \max_{k=1,..,k_p} \{ z_{k+1/2} \, z_{p-k}, \,  z_{k} \,
         z_{p-k+1/2} \}, \quad z_p := {m_p \over \Gamma(a\,p+1/2)}. \]
We may then rewrite~(\ref{mp6}) as
      \beqn \label{mp7}
      \qquad {d z_p \over dt} \le \alpha(\EE) \, \left(A \, \gamma_p  \,
{\Gamma (a\,p+a/2 + 1)
      \over  \Gamma(ap+1/2)} \, Z_p
      - (1-\gamma_p) \, \Gamma(a\,p+1/2)^{1/2p} \, z_p^{1+1/2p} \right)
      \eeqn
for any $p=3/2,2,\dots$ On the one hand, from (\ref{mp4}), there exists
$A'$ such that
      \beqn \label{mp8}
      A \, \gamma_p \, {\Gamma(ap+a/2+1) \over  \Gamma(ap+1/2)}  \le A'
\,
p^{a/2-1/2} \qquad
      \forall \, p=3/2,2, \dots
      \eeqn
On the other hand, thanks to Stirling's formula $n! \sim n^n \,
e^{-n} \, \sqrt{2\pi n}$ when $n\to\infty$ and the estimate (\ref{mp4}),
there exists $A'' > 0$ such that
      \beqn \label{mp9}
      (1-\gamma_p) \, \Gamma(a\,p+1/2)^{1/2p} \ge A'' \, p^{a/2}  \qquad
\forall
      \,  p = 3/2,2, \dots
      \eeqn
Gathering~(\ref{mp7}), (\ref{mp8}) and~(\ref{mp9}), we obtain the
differential inequality
      \beqn \label{mp10}
      \qquad {d z_p \over dt} \le \alpha(\EE) \, \left( A'  \,
p^{a/2-1/2} \, Z_p -
A'' \, p^{a/2} \, z_p^{1+1/2p} \right)
    \eeqn
for any $p = 3/2,2, \dots$

\medskip\noindent{\sl Step 2. Proof of (\ref{MomentUnif2}). }
On the one hand, we remark, by an induction argument,  that taking $p_0
:= \max\{3/2,(2A'/A'')^2\}$, the
sequence of functions $z_p := x^p$ is a sequence of  supersolutions
of~(\ref{mp10}) for any $x > 0$ and
for $p \ge p_0$. On the other hand, choosing $x_0$ large enough, which
may depend on $p_0$, with have from (i)  that  the sequence of
functions $z_p := x^p$ is a
sequence of  supersolutions of~(\ref{mp10}) for any $x \ge x_0$ and
for $p \in \{ 3/2, \dots, p_0 \}$. As a consequence, since $z_p$ for $p
= 0,1/2,1$ are bounded by
$\| f_{\mbox{\scriptsize{in}}}\|_{L^1_2}$, we have proved that there
exists $x_0$ such that
the set
    \beqn
    \label{superS}
    \Cc_x := \left\{ z = (z_p); \quad z_p \le x^p \,\, \, \forall \, p
\in {1
    \over 2} \, \NN \right\}
    \eeqn
is invariant under the flow generated by the Boltzmann equation
for any $x \ge x_0$: if $f(t_1) \in \Cc_x$ then $f(t_2) \in \Cc_x$ for
any $t_2 \ge t_1$.
We set $a := 2/\eta \ge 1$. Noticing that
     \beqn
     \int_{\R^N} f(v) \, e^{r\,|v|^\eta} \, dv = \sum_{k=0}^\infty {r^k
\over k!} \, m_{k \, \eta/2}
     \eeqn
     we get, from the assumption made on $f_{\mbox{\scriptsize{in}}}$,
that
     $$
     m_{k/a}(0) \le C_0 \, {k! \over r^k} \quad \forall \, k \in \N.
     $$
Since we may assume $r \in (0,1]$, the function $y \mapsto C_0 \,
\Gamma(y+1) r^{-y}$
is increasing, and we deduce by H\"older's
inequality that for any $p$
     $$
     m_p (0) \le C_0 \, {\ell_p! \over r^{\ell_p}} \le C_0 \,
{\Gamma(ap+2)
     \over r^{ap+2}} \quad\hbox{with}\quad \ell_p := [a \, p]+ 1.
     $$
 From the definition of $z_p$ we deduce
    \beqn
    \label{xp0lex0p}
    z_p(0) \le C_0\, {ap \, (ap+1) \over r^{ap+2}}\le x_1^p
    \eeqn
for any $p$ and for some constant $x_1 \in (0,\infty)$. Choosing $x
:= \max\{x_0,x_1\}$ we get  from (\ref{superS}) and (\ref{xp0lex0p}) that
$z_p(t) \le x^p$ $\forall \, t \in [0,T_c)$  for any $p$.Therefore, we
have
    $$
    m_p(t) \le \Gamma(ap+1/2) \, x^p \qquad \forall \, p = 3/2, 2, \dots,
\,\,\, \forall \, t \in [0,T_c).
    $$
The function $y \mapsto \Gamma(y+1/2) \, x^y$ being increasing, we
deduce from H\"older's inequality that for any $k \in \N^*$ that
$m_{k/a}(t) \le \Gamma(ap+1/2) \, x^p \le \Gamma(k+a/2+1/2) \,
x^{k/a+1/2}$
with $p := [2k/a]/2+ 1/2$.
For $r' < 2 x^{-1/a} (1+a)^{-1}$ we conclude
    \bean
    \forall \, t \in [0,T_c) \quad
    \int_{\R^N}f(t,v) \, e^{r' \, |v|^\eta} \, dv \le  \sum_{k=0}^\infty
    {\Gamma(k+a/2+1/2) \over k!} \, {x^{k/a+1/2}}
    \, (r')^k \\
    \le C \, \sum_{k=0}^\infty \left( \left(\frac{a+1}2\right) x^{1/a} r' \right)^k < +\infty
    \eean
from which (\ref{MomentUnif2}) follows.

\medskip\noindent{\sl Step 3.  Proof of (\ref{expeta}). } 
Let us fix $\tau \in (0,T_c)$. We claim that there exists $x$ large enough and some
increasing  sequence of times $(t_p)_{p \ge p_0}$ which are bounded by
$\tau$ such that for any $p$
     \beqn\label{wple1}
\forall \, t \in [t_p,T_c) \quad  z_p (t) \le x^p.
\eeqn
We already know by classical arguments (see \cite{MW99,CVhand}) that for
$p_0$ (defined at the beginning of Step 2) there exists $x_1$, larger
than $x_0$ defined in (\ref{superS}), such that (\ref{wple1}) holds for
any $p \le p_0$ and $t_p = \tau/2$.
We then argue by induction, assuming that for $p \ge p_0$ there holds:
    \bear\label{wkle1}
    & & z_k \le x^k \quad \hbox{on}\quad [t_{p-1/2},T_c) \quad \forall \,
k
    \le p-1/2 \\
    \label{wpge1}
    &&  z_p \ge x^p \quad \hbox{on}\quad [t_{p-1/2},t_p),
    \eear
for some $x \ge x_1$ to be defined.
If (\ref{wpge1}) does not hold, there is nothing to prove thanks to
Step 2. Gathering (\ref{wkle1}), (\ref{wpge1}) with (\ref{mp10}) we get
from the definition of $p_0$ and the fact that $\EE(t) \in
[\EE(\tau),\EE(0)]$ so that $\alpha (\EE) \ge \alpha_0 > 0$
     \beqn \label{mpstep41}
     \qquad {d z_p \over dt} \le - \alpha_0 \,  {A'' \over 2}  \, p^{a/2}
     \, z_p^{1+1/2p}  \quad\hbox{on}\quad (t_{p-1/2},t_p).
     \eeqn
Integrating this differential inequality we obtain
    $$
    -  z_p^{-{1\over 2p}} (t_p) \le z_p^{-{1\over 2p}} (t_{p-1/2}) -
    z_p^{-{1\over 2p}} (t_p) \le
    -{1 \over 2p}\, {A'' \, \alpha_0 \over 2}  \, p^{a/2} \, (t_p -
    t_{p-1/2}).
    $$
Defining $(t_p)$ in the following way:
    $$
    t_0 := {\tau \over 2}, \quad t_p := t_{p-1/2} + {\tau \over 2} \, {
    p^{1-a/2} \over s_a},
    \quad s_a := \sum_{p=0}^\infty p^{1-a/2}
    $$
and defining $x_2 := (8 \, s_a)^2/(A'' \, \alpha_0 \, \tau)^2$ we have
then proved $z_p(t_p) \le x_2^p$ and therefore $z_p (t) \le x^p$ for any
$t  \ge (t_p,T_c)$ with $x = \max \{ x_1,x_2 \}$  thanks to Step 2.
Setting $a:= 2/\eta > 4$ ($\eta < 1/2$) we have
\beqn\label{defTeta}
    \sum_{k=0}^\infty t_{1+k/2} \le \tau
\eeqn
and we conclude as in the end of Step 2.
\qed

\subsection{Stability estimate in $L^1 _2$ and proof of the uniqueness
part of Theorem~\ref{Lcase}}

      \begin{prop}\label{uniqueness}
      Assume that $B$ satisfies~(\ref{shapeB1t})--(\ref{shapeB3t}). For any
      two solutions $f$ and $g$ of the inelastic
      Boltzmann equation~(\ref{eqB1})-(\ref{eqB2}) on $[0,T]$ ($T \le
T_*$) we have
        \beqn \label{eq:stab}
        \qquad\,\,\, {d \over dt} \int_{\RR^N} \! |f-g| \, (1+|v|^2) \,
dv
        \le C \, \gamma_* \! \int_{\RR^N} \! (f+g) \, (1+|v|^3) \, dv
        \int_{\RR^N} |f-g| \, (1+|v|^2) \, dv.
        \eeqn
      We deduce that there is $C _T >0$ depending on $B$ and
      $\sup_{t\in[0,T]} \|f+g\|_{L^1 _3}$ such that
        \[ \forall \, t \in [0,T], \ \ \
           \|f_t -g_t \|_{L^1 _2} \le  \|f_{\mbox{\scriptsize{{\em
in}}}} -g_{\mbox{\scriptsize{{\em in}}}} \|_{L^1 _2} \,
           e^{C_T t}. \]
      In particular, there exists at most one solution to the Cauchy
      problem for the inelastic Boltzmann equation in
      $C([0,T];L^1_2) \cap L^1(0,T;L^1_3)$.
      \end{prop}

\smallskip\noindent
{\sl Proof of Proposition~\ref{uniqueness}. }
We multiply the equation satisfied by $(f-g)$ by
$\phi(t,y) = \hbox{sgn}(f(t,y)-g(t,y)) \, k$ , where $k=(1+|v|^2)$.
Using the chain rule~(\ref{renormBeq}), we get for all $t \geq 0$
      \bean
      \frac{d}{dt} \int_{\RR^N} |f-g| \,  k  \, dv
      &=& \frac{1}{2} \, \int_{\RR^N \times \RR^N \times D} \, \left[
      (f-g) g_* + f
      (f_* - g_*)\right] \\
      && \hspace{2cm} (\phi' + \phi'_* - \phi -
      \phi_*)\, B(t,u;dz) \,  dv_* \, dv \\
      &=& \frac{1}{2} \, \int_{\RR^N \times \RR^N \times D}
         (f-g) \, (f_* + g_*) \\
      && \hspace{2cm} (\phi' + \phi'_* - \phi -
      \phi_*) \, B(t,u;dz) \, dv_* \, dv \\
      &\leq& \frac{1}{2} \, \int_{\RR^N \times \RR^N \times D} |f-g| \,
      (f_* + g_*) \\
      && \hspace{2cm} (k' + k'_* - k + k_*) \, B(t,u;dz) \, dv_* \, dv,
      \eean
where we have just use the symmetry hypothesis~(\ref{shapeB1t}),
(\ref{shapeB2t}) on $B$
and a change of variable $(v,v_*) \to (v_*,v)$. Then, thanks to the
bounds
(\ref{shapeB1t}), (\ref{shapeB3t})  we deduce
      \bean
      \frac{d}{dt} \int_{\RR^N} |f-g| \, k \, dv
      &\leq& \gamma_* \, \int_{\RR^N \times \RR^N}  |u| \, |f-g| \,
      (f_* + g_*) \, k_* \,  dv_*dv \\
      &\leq& \gamma_* \, \int_{\RR^N} |f-g| \, k \, dv \int_{\RR^N}  (f_*
+ g_*) \,
      k_*^{3/2} \,  dv_*
      \eean
which yields the differential inequality~(\ref{eq:stab}).
The end of the proof is straightforward by a Gronwall argument. \qed

\medskip
The uniqueness in $C([0,T);L^1_2) \cap L^1(0,T;L^1_3)$ as stated in
Theorem \ref{Lcase} is given by Proposition~\ref{uniqueness}.

\subsection{Sketch of the proof of the existence  part of
Theorem~\ref{Lcase}}\label{subsec:lin}

As for the existence part, we briefly sketch the proof. We follow a
method introduced in~\cite{MW99}
and developed in~\cite{FM**}. We split the proof into three steps.

\medskip\noindent
{\sl Step 1. }
Let us first consider an initial datum $f_{\mbox{\scriptsize{in}}}$
satisfying
(\ref{initialcond}) with $q = 4$ and let us define the truncated
collision rates
$B_n = B \, {\bf 1}_{|u| \le n}$. The associated collision operators
$Q_n$ are
bounded in any $L^1_q$, $q \ge 1$, and are Lipschitz in $L^1_2$ on any
bounded subset of $L^1_2$. Therefore following a classical argument
from Arkeryd, see \cite{Ar72}, we can use  the Banach fixed point
Theorem and obtain the existence of a solution $0 \le f_n \in
C([0,T];L^1_2) \cap L^\infty(0,T;L^1_4)$ for any $T > 0$, to the
associated Boltzmann equation (\ref{eqB1})-(\ref{eqB2}), which satisfies
(\ref{consmq})-(\ref{decE}).

\medskip\noindent
{\sl Step 2. } From Proposition~\ref{propboundL3}, for any $T > 0$,
there exists $C_T$ such that
\[
\sup_{[0,T]}\|f_n \|_{L^1_4} \le C_T.
\]
Moreover, coming back to the proof of Proposition~\ref{uniqueness} (see
also the first step in the proof of \cite[Theorem 2.6]{FM**}), we may
establish the differential inequality
\[
{d \over dt}\|f_n - f_m \|_{L^1_2} \le C_1 \, \|f_n + f_m \|_{L^1_3}
\, \|f_n - f_m \|_{L^1_2}
+ {C_2 \over n}\,   \|f_n + f_m \|_{L^1_4}^2
\]
for any integers $m \ge n$. Gathering these two informations we easily
deduce that $(f_n)$ is a Cauchy sequence in $C([0,T];L^1_2)$ for any $T
  > 0$.
Denoting by $f \in C([0,T];L^1_2) \cap L^\infty(0,T;L^1_4)$
its limit, we obtain that $f$ is a solution to the Boltzmann equation
(\ref{eqB1})-(\ref{eqB2}) associated to the collision rate $B$ and the
initial datum $f_{\mbox{\scriptsize{in}}}$ by passing to the limit in
the weak formulation
(\ref{dualBeq}) of the Boltzmann equation written for $f_n$.

\medskip\noindent
{\sl Step 3. } When the initial datum $f_{\mbox{\scriptsize{in}}}$
satisfies
(\ref{initialcond}) with $q = 3$ we introduce the sequence of initial
data $f_{\mbox{{\scriptsize in}},\ell} := f_{\mbox{\scriptsize{in}}}\, {\bf 1}_{\{|v|\le
\ell \} }$. Since
$f_{\mbox{{\scriptsize in}},\ell} \in L^1_4$, the preceding  step give the existence of a
sequence of solutions $f_\ell \in C([0,T];L^1_2) \cap
L^\infty(0,T;L^1_3)$ for any $T > 0$ to the Boltzmann
equation~(\ref{eqB1})-(\ref{eqB2})
associated to the  initial datum $f_{\mbox{{\scriptsize in}},\ell}$. From
Proposition~\ref{propboundL3}, for any $T > 0$, there exists $C_T$
such that
\[
\sup_{[0,T]}\|f_\ell \|_{L^1_3} \le C_T.
\]
Thanks to (\ref{eq:stab}) we establish that $(f_\ell)$ is a Cauchy
sequence in $C([0,T];L^1_2)$
and we conclude as before. \qed

\begin{rem} Note here that an alternative path to the proof of existence
could have been the use of the result of propagation of Orlicz norm
which shows that the solution is uniformly bounded
for $t \in [0,T]$ in a certain Orlicz space. Together with
the propagation of moments and Dunford-Pettis Lemma, it would yield
the existence of a solution by classical approximation arguments and
weak stability
results as presented below.
More generally the propagation of Orlicz norm by the collision
operator can be seen as a new tool (as well as a clarification)
for the theory of solutions to the spatially homogeneous Boltzmann
equation with no entropy bound, as in the inelastic case, or
in the elastic case when the initial datum has infinite  entropy, see
also \cite{Ar72,MW99} where
other strategies of proof are presented.
\end{rem}


\section{Proof of the Cauchy theorem for coupled collision
rate}\label{sec:cauchy:NL}
\setcounter{equation}{0}
\setcounter{theo}{0}


\subsection{Weak stability and proof of the existence part of Theorem
\ref{NLcase}} \label{subsec:NL}

      \begin{prop}\label{weakstab}
      Consider a sequence $B_n = B_n(t,u;dz)$ of collision rates
      satisfying the structure
conditions~(\ref{shapeB1t})-(\ref{shapeB2t})
      and the uniform bound
        \[  
        0 \le \gamma_n(t) \le \gamma_T \quad \forall \, t \in [0,T],\quad
        \forall \, n \in \N^*,
        \]
      and let us denote by $f_n \in C([0,T);L^1_2) \cap
L^\infty(0,T;L^1_3)$ the
      solution associated to $B_n$ thanks to the existence result of the
preceding
      section (existence and uniqueness part of Theorem~\ref{Lcase} and
Remark~\ref{remLinear}
      4th point). Assume furthermore that $(f_n)$ belongs to a weak
compact set
      of $L^1((0,T) \times \R^N)$ and that   there exists a collision
rate
      $B$ satisfying (\ref{shapeB1t})-(\ref{shapeB2t})-(\ref{shapeB3t}) and such that for
any
      $\psi \in C_c(\RR^N)$
      \[
      \gamma_n \ \to \ \gamma \quad\hbox{and}\quad
      \int_D \psi(v')  \, b_n (t,u;dz) \ \to \ \int_D \psi(v')  \,
      b(t,u;dz) \quad a.e.
      \]
Then there exists a function $f \in C([0,T);L^1_2)\cap
L^\infty(0,T;L^1_3)$ and a subsequence $f_{n_k}$  such that
\[ f_{n_k} \, \wto \, f \quad\hbox{weakly in}\quad L^1((0,T) \times
\R^N), \]
and $f$ is a solution to the Boltzmann equation
(\ref{eqB1})-(\ref{eqB2}) associated to $B$.
\end{prop}

\smallskip
Such a stability/compactness result is very classical and we refer to
\cite{Ar72,DPL89}
for  its proof.

\medskip\noindent
{\sl Proof of the existence part of Theorem \ref{NLcase}. }
We assume without restriction that there exists a
decreasing function $\alpha_0$ such that $\alpha \le \alpha_0$ on
$[0,\EE_{\mbox{\scriptsize{in}}}]$. We proceed in three steps.

\smallskip\noindent{\sl Step 1. }
We start with some {\it a priori} bounds. We set $Y_3 := \| f
\|_{L^1_3}$.
  From the Povzner inequality
(\ref{Y3gal}) (with $\gamma(t) = \alpha(\EE(t)))$ and the dissipation of
energy equation (\ref{eqdiffEE}), we have
      \beqn
      \label{IDy3}
      {d \over dt} Y_3 \le C_1 \, \alpha_0(\EE) \, Y_3, \quad Y_3(0)      
Y_3(f_{\mbox{\scriptsize{in}}})
      \eeqn
and
      \beqn
      \label{IDEE}
      {d \over dt} \EE \ge -  C_1 \, \alpha_0(\EE) \, Y_3,  \quad \EE(0)
=\EE_{\mbox{\scriptsize{in}}},
      \eeqn
for some constant $C_1$ (which depends on
$\EE_{\mbox{\scriptsize{in}}}$).
There exists $T_*$ such that any solution $(Y_3,\EE)$ to the above
differential inequalities system is defined on $[0,T_*]$ and satisfies
      \beqn
      \label{bornesY3&EE}
      \sup_{[0,T_*]} Y_3(t) \le 2 \, Y_3(f_{\mbox{\scriptsize{in}}}),
\quad \inf_{[0,T_*]}
      \EE(t) \ge \EE_{\mbox{\scriptsize{in}}}/2.
      \eeqn
More precisely, we choose $T_*$ such that
      \[
      \label{defT*}
      C_1 \, \alpha_0(\EE_{\mbox{\scriptsize{in}}}/2) \, T_* \le
Y_3(f_{\mbox{\scriptsize{in}}})
      \quad\hbox{and}\quad
      C_1 \, \alpha_0(\EE_{\mbox{\scriptsize{in}}}/2) 2 \,
Y_3(f_{\mbox{\scriptsize{in}}}) \, T_* \le
\EE_{\mbox{\scriptsize{in}}}/2,
      \]
in such a way that if $(Y_3,\EE)$ satisfies $Y_3 \le 2 \,
Y_3(f_{\mbox{\scriptsize{in}}})$
and
(\ref{IDEE}) on $(0,T_*)$ or if $(Y_3,\EE)$ satisfies $\EE \ge
\EE_{\mbox{\scriptsize{in}}}/2$ and (\ref{IDy3}) on
$(0,T_*)$ then (\ref{bornesY3&EE}) holds. We introduce
      $$
      X := \Big\{ \EE \in C([0,T_*]), \,\,
\EE_{\mbox{\scriptsize{in}}}/2 \le \EE(t) \le
\EE_{\mbox{\scriptsize{in}}}
      \, \, \hbox{ on } \,\, (0,T_*) \Big\}.
      $$

\smallskip\noindent
{\sl  Step 2.}
Let us consider a function $\EE_1 \in X$ and define $B_2(t,u;dz):=B(\EE_1(t),u;dz)$.
From assumptions (\ref{shapeB1})-(\ref{shapeB2})-(\ref{shapeB3})-(\ref{shapeB4} 
we may write
      $$
      B_2(t,u;dz) = |u| \, \gamma_2(t) \, b_2(t,u;dz)
      $$
where $b_2$ is a probability measure and $\gamma_2(t)$ satisfies
      $$
      \gamma_2(t) = \alpha(\EE_1(t)) \le
\alpha_0(\EE_{\mbox{\scriptsize{in}}}/2) < + \infty
      \qquad \forall \, t \in [0,T_*].
      $$
Thanks to Theorem \ref{Lcase} there exists a unique solution
$f_2 \in C([0,T_*];L^1_2) \cap L^\infty(0,T_*;L^1_3)$ to the
Boltzmann equation (\ref{eqB1})-(\ref{eqB2}) associated to the
collision rate $B_2$ and we set $\EE_2 := \EE(f_2)$.
In such a way we have defined
a map $\Phi : X \to X$, $\Phi(\EE_1) = \EE_2$.

In order to apply the Schauder fixed point Theorem, we aim
to prove that $\Phi$ is continuous and compact from $X$ to $X$.
Consider $(\EE^n_1)$ a sequence of $X$ which
uniformly converges to $\EE_1$. Since $(\EE^n _1)$ belongs
to the compact set
$[\EE_{\mbox{\scriptsize{in}}}/2,\EE_{\mbox{\scriptsize{in}}}]$ for any
$n$ and any $t \in [0,T_*]$, we deduce by applying
Corollary~\ref{cor:Orl} to
the sequence $(f^n_2)$ associated to $B_2^n(t,u;dz) B(\EE^n_1(t),u;dz)$ that
     \beqn\label{Lambdaf2n}
     \forall \, n \ge 0, \quad \sup_{[0,T_*]}\int_{\R^N}
\Lambda(f_2^n(t,v)) \, dv \le C_2,
     \eeqn
for a superlinear function $\Lambda$ and a constant $C_2>0$.
Moreover, from Proposition~\ref{propboundL3} we have
\beqn\label{L3f2n}
\forall \, n \ge 0, \quad \sup_{[0,T_*]}\int_{\R^N} f_2^n(t,v) \, |v|^3
\, dv \le C_3
\eeqn
for some constant $C_3>0$.

On the one hand, gathering (\ref{Lambdaf2n}), (\ref{L3f2n}) and using
the Dunford-Pettis
Lemma, we obtain that $(f_2^n)$ belongs to a weak compact set of
$L^1((0,T_*) \times \R^3)$. Propositon~\ref{weakstab} then implies that
there exists $f_2 \in C([0,T_*];L^1_2) \cap
L^\infty(0,T_*;L^1_3)$ such that, up to a subsequence, $f_2^n \wto f_2$
weakly in
$L^1(0,T;L^1_2)$ and $f_2$ is a solution to the Boltzmann equation
associated to $B_2(t,u;dz) = B(\EE_1(t),u;dz)$. Since this limit is
unique by
the previous study, the whole sequence $(f_2^n)$ converges weakly to
$f_2$, and
in particular
    \beqn\label{En2wtoE2}
    \EE^n_2 \,\, \wto \,\, \EE_2 \quad\hbox{weakly in}\quad L^1(0,T)
    \eeqn
where $\EE_2$ is the kinetic energy of $f_2$.

On the other hand, there holds
    \[
    {d \over dt} \EE^n_2 = -  \, \int_{\RR^N \times \RR^N}
    f_2^n \, f_{2*} ^n\, |u|^3 \, \Delta(\EE_1^n,u) \, dvdv_*
    =: - D^n_2.
    \]
Since $\Delta(\EE_1^n,u) \le \alpha(\EE_1^n)/4  \le
\alpha_0(\EE_{\mbox{\scriptsize{in}}}/2)/4$,
we deduce from (\ref{propboundL3}) that $D^n_2$ is bounded in
$L^\infty(0,T)$ which in turn implies
    \beqn\label{En2Lip}
    \|\EE^n_2 \|_{W^{1,\infty}(0,T)}\le C_4.
    \eeqn
  From Ascoli's Theorem we infer that the sequence $(\EE^n_2)$ belongs
to a compact set
of $C([0,T])$. Since the cluster points for the uniform norm are
included in the set of cluster points for the $L^1(0,T)$ weak topology, it then
follows from
(\ref{En2wtoE2}) that $\Phi(\EE^n_1) = \EE(f^n_2)$ converges to
$\Phi(\EE_1) = \EE(f_2)$ for the uniform norm on $C([0,T])$, which
ends the proof of the continuity of $\Phi$. Of course, the {\it a
priori}
bound (\ref{En2Lip}) and Ascoli's Theorem also imply that $\Phi$ is a
compact map on $X$. We may thus use
the Schauder fixed point Theorem to conclude to the
existence of at least one $\bar\EE \in X$ such that $\Phi(\bar\EE) = \bar\EE$. Then,
the solution $\bar f \in C([0,T_*];L^1_2)
\cap L^\infty(0,T_*;L^1_3)$ to the Boltzmann equation associated to
$\bar B(t,u;dz) := B(\bar\EE(t),u;dz)$ satisfies
$$
     \int_{\RR^N} \bar f(t,v) \, |v|^2 \, dv = \Phi(\bar \EE)(t) = \bar
\EE(t)
$$
and therefore $\bar f$ is a solution to the Boltzmann
equation associated to $B$ in \linebreak $C([0,T_*];L^1_2) \cap
L^\infty(0,T_*;L^1_3)$.

\smallskip\noindent{\sl Step 3. }
We then consider the class of solution $f : (0,T_1) \to L^1_3$ such that
$f \in C([0,T];L^1_2) \cap L^\infty(0,T;L^1_3)$ for any $T \in
(0,T_1)$, $\EE$ is decreasing, $f$
is mass conserving.
By Zorn's Lemma, there exists a maximal interval $[0,T_c)$ such that
      $$
      (T_c < \infty \hbox{ and } \EE(t) \to 0 \hbox{ when } t \to T_c)
      \quad\hbox{ or }\quad
      T_c = + \infty.
      $$
      In order to end the proof, the only thing one has to remark is that
if $T_c < + \infty$ and
$\displaystyle{\lim_{t \nearrow T_c} \EE(t) = \EE_c > 0}$,
then $\displaystyle{\lim_{t \nearrow T_c} Y_3(t) < \infty}$ (by
(\ref{IDy3})) so that
$f \in C([0,T_c];L^1_2) \cap L^\infty(0,T_c;L^1_3)$ and we
may extends the solution $f$ to a larger time interval. \qed

\subsection{Strong stability and uniqueness part of
Theorem~\ref{NLcase}}

In this subsection we give a quantitative stability
result in strong sense, under the additional assumption of some
smoothness on the
initial datum and the collision rate.
Let us first prove a simple result of propagation of the total variation
of the distribution.

     \begin{prop} \label{sobolev}
     Let $B$ be a collision rate satisfying assumptions
     (\ref{shapeB1t})-(\ref{shapeB2t})-(\ref{shapeB3t})
     and $0 \le f_{\mbox{\scriptsize{{\em in}}}} \in BV_4 \cap L^1_5$ an
initial datum.
     Then there exists $C_{T_*}$, depending
     on $\gamma_*$ and $\|f_{\mbox{\scriptsize{{\em in}}}}\|_{L^1 _5}$,
such that
     any solution  $f \in C([0,T_*],L^1 _2) \cap L^\infty(0,T_*,L^1 _3)$
to
     the Boltzmann equation  constructed in the previous step satisfies
        \[ \forall \, t \in [0,T_*], \ \ \
           \|f_t\|_{BV_4} \le \|f_{\mbox{\scriptsize{{\em
in}}}}\|_{BV_4} \, e^{C _{T_*}
           \, t}. \]
      \end{prop}

\smallskip \noindent
{\sl Proof of Proposition~\ref{sobolev}.}
The proof is based on the same kind of Povzner inequality as above.
Let us first prove the estimate by {\em a priori} approach, for the
sake of clearness.
We have the following formula for the differential of $Q$:
      \[ \nabla_v Q(f,f) = Q(\nabla_v f,f) + Q(f, \nabla_v f ). \]
This property is proved in the elastic case in~\cite{CVhand} but it
is strictly related to the invariance property of the collision operator
      \[ \tau_h Q(f,f) = Q(\tau_hf, \tau_h f) \]
where the translation operator $\tau_h$ is defined by
      \[ \forall \ v \ \in \ \R^N, \ \ \ \tau_h f (v) = f(v-h). \]
It is easily seen that it remains true in the
inelastic case under our assumptions. The propagation of the
$L^1 _5$ norm has already been established.
Then we estimate the time derivative of the $L^1 _4$ norm of the
gradient along the flow:
      \bean
      \frac{d}{dt} \|\nabla_v f_t  \|_{L^1_4}
      &=&
         \int_{\RR^N \times \RR^N \times D} f \, (\nabla_v f_*)
         \Big[ (1+|v'|^4) \, \mbox{sgn}(\nabla_v f)' + (1+|v'_*|^4) \,
          \mbox{sgn}(\nabla_v f)' _* \\
       &&         - (1+|v|^4) \, \mbox{sgn}(\nabla_v f) - (1+|v_*|^4) \,
       \mbox{sgn}(\nabla_v f)_*
         \Big] \, B \, dv \, dv_* \\
      &\le&
         \int_{\RR^N \times \RR^N \times D} f \,  |\nabla_v f_*|
         \Big[ (1+|v'|^4) + (1+|v'_*|^4) - (1+|v|^4)  \\
         &&- (1+|v_*|^4) \Big] \, B \, dv \, dv_*
         + 4 \, \gamma_* \, \|f_t (1+|v|^5) \|_{L^1} \, \| \nabla_v f
         (1+|v|) \|_{L^1} \\
         &\le& C \, \|f_t\|_{L^1_5} \, \| \nabla_v f \|_{L^1 _4}
      \eean
using a Povzner inequality as in (\ref{povzner1}).
This shows the {\it a priori} propagation of the $BV_4$ norm by a
Gronwall argument.

Now let us explain how to obtain the same estimate by {\em a posteriori}
approach. First concerning the {\it a posteriori} propagation of the
$L^1 _5$
norm, it is similar to the method in~\cite{MW99} and does not lead to
any difficulty. Concerning the propagation of $BV_4$ norm,
we look at some ``discretized derivative''.
Let us denote $k=\mbox{sgn}(\tau_h f -f) \, (1+|v|^4)$.
We can compute by the chain rule the following time derivative
(using the invariance property of the collision operator)
      \bean
      \frac{d}{dt} \|\tau_h f_t - f_t \|_{L^1_4}
      &=&
         \int_{\RR^N \times \RR^N \times D} \big( \tau_h f \tau_h f_* - f
f_*) \,
         \big[ k' - k \big] \, B \, dv \, dv_* \\
      &=&
         \int_{\RR^N \times \RR^N \times D} (\tau_h f -f) f_* \,
         \big[ k' + k' _* - k - k_* \big] \, B \, dv \, dv_* \\
      &&  + \frac{1}{2} \, \int_{\RR^N \times \RR^N \times D} (\tau_h f
-f) (\tau_h f_* - f_*) \,
         \big[ k' + k' _* - k - k_* \big] \, B \, dv \, dv_* \\
      &\le&  \int_{\RR^N \times \RR^N \times D} |\tau_h f -f| f_* \,
         \Big[ |v'|^4 + |v' _*|^4 - |v|^4 + |v_*|^4 \Big] \, B \, dv \,
dv_* \\
      && + \frac{1}{2} \, \int_{\RR^N \times \RR^N \times D} |\tau_h f
-f|
|\tau_h f_* - f_*| \\
      && \hspace{3cm} \Big[ |v'|^4 + |v' _*|^4 + |v|^4 + |v_*|^4 \Big] \,
B \, dv \, dv_*.
      \eean
Then using the same rough Povzner inequality as in the proof of
Proposition~\ref{propboundL3},
we have
\[ \Big[ |v'|^4+ |v' _*|^4 + |v|^4 + |v_*|^4 \Big] |v-v_*| \le C \, \Big[
             (1+|v|^4)(1+|v_*|^5) + (1+|v_*|^4)(1+|v|^5) \Big]. \]
Hence we deduce that
     \[ \frac{d}{dt} \|\tau_h f_t - f_t \|_{L^1_4} \le C \, \gamma_* \,
\|\tau_h f_t - f_t \|_{L^1_4}
                          \left[ \|f\|_{L^1 _5} + \|\tau_h f_t - f_t
\|_{L^1_5} \right] \]
and for $|h| \le 1$, we deduce
     \[ \frac{d}{dt} \|\tau_h f_t - f_t \|_{L^1_4} \le C \, \gamma_* \,
\|\tau_h f_t - f_t \|_{L^1_4}
                          \|f\|_{L^1 _5}. \]
By a Gronwall argument it shows for any $|h| \le 1$ that
      \[ \forall \, t \in [0,T_*], \ \ \
         \|\tau_h f_t - f_t\|_{L^1 _4} \le \|\tau_h
f_{\mbox{\scriptsize{in}}} - f_{\mbox{\scriptsize{in}}}\|_{L^1
_4} \, e^{C _{T_*} \, t} \]
for a constant $C_{T_*}$ depending on $\gamma_*$
and $\sup_{t \in [0,T_*]} \|f_t\|_{L^1 _5}$.
By dividing by $h$ and letting $h$ goes to $0$, we conclude that
     \[ \forall \, t \in [0,T_*], \ \ \
         \|\nabla_v f_t\|_{M^1 _4} \le \|\nabla_v
f_{\mbox{\scriptsize{in}}} \|_{M^1 _4} \,
e^{C _{T_*} \, t} \]
which ends the proof.
\qed

\medskip
Assume now that the collision rate satisfies
(\ref{shapeB1})-(\ref{shapeB2})-(\ref{shapeB3})-(\ref{shapeB4})
plus the additional assumption {\bf H1}. 
Let us take $f_{\mbox{\scriptsize{in}}} \in BV_4 \cap L^1_5$ and let us consider
two solutions $f,g \in C([0,T_c];L^1_2) \cap L^\infty(0,T;L^1_3)$
constructed by the
previous steps. For these two solutions the function $e(\EE)$ is locally
Lipschitz, so is the function $\Delta(\EE)$ and
the differential equation~(\ref{eqdiffEE}) satisfied by $\EE(f_t)$ on
$[0,T_*]$ implies that it is bounded from below on this interval.
Thus thanks to the continuity of $\alpha$, the assumptions of
Proposition~\ref{sobolev} are satisfied, and thus the
$BV_4$ norm is bounded on any time interval $[0,T_*] \subset [0,T_c)$ .
       
\begin{prop} \label{uniqNL}
Let $B$ be a collision rate satisfying (\ref{shapeB1})-(\ref{shapeB2})-(\ref{shapeB3})-(\ref{shapeB4}) 
plus the additionnal assumption {\bf H1}.
Let $f,g \in C([0,T_*];L^1_2) \cap L^\infty(0,T_*;L^1_3)$ be two
solutions with mass $1$ and momentum $0$, with initial data $f_{\mbox{\scriptsize{{\em in}}}}$ 
and $g_{\mbox{\scriptsize{{\em in}}}}$, and such that $\EE(f(t,.))$, $\EE(g(t,.)) \in K$ on $[0,T_*]$
     with $K$ a compact of $(0,+\infty)$ and
\[ 
\forall \, t \in [0,T_*], \ \ \ \|f(t,.)\|_{BV_4}, \|g(t,.)\|_{BV_4} \le C_{T_*}. 
\]
Then there is a constant $C' _{T_*}$ depending on $B$, $K$ and $C_{T_*}$ such that
\[ 
\forall \, t \in [0,T_*], \ \ \ \|f(t,.) - g(t,.)\|_{L^1 _2} 
 \le \|f_{\mbox{\scriptsize{{\em in}}}} - g_{\mbox{\scriptsize{{\em in}}}}\|_{L^1 _2} \, e^{C' _{T_*} \, t}. 
\]
\end{prop}

We need the following geometrical lemma which is a more accurate version of 
Lemma~\ref{utow} when the collision process is of the 
generalized visco-elastic type (\ref{defCue},\ref{bCue}). 

\begin{lem}\label{v&v*tovprim} 
For any $e \in (0,1]$ and $\sigma \in \Sph^{N-1}$ we define 
\bear\label{v*tovprim}
&&\phi^*_e = \phi^*_{e,v,\sigma} : \R^N \to \R^N, \quad v_* \mapsto  
v' =  v + {1 + e \over 4} \, \Phi_{\sigma}(v_*-v) \\ \label{vtovprim}
&&\phi_e = \phi_{e,v_*,\sigma} : \R^N \to \R^N, \quad v \mapsto  
v' =  v_* + {3-e \over 4} \, \Phi_{r_e \sigma}(v-v_*), \,\, r_e = {1 + e \over 3-e},
\eear
(where $\Phi_z$ was defined in Lemma~\ref{utow}) 
and the Jacobian functions $J^*_e =\mbox{{\em det}} \, (D \,  \phi^*_{e,v,\sigma})$, 
$J_e = \mbox{{\em det}} \, (D \,  \phi_{e,v_*,\sigma})$.

Then for any $\gamma \in (-1,1)$, $\phi^*_e$ defines a $C^\infty$-diffeomorphism 
from $v + \Omega_\gamma$ onto $v + \Omega_{\omega^*(\gamma)}$ 
with $\omega^*(\gamma) = ((1+\gamma)/2)^{1/2}$ and $\phi_e$ 
defines a $C^\infty$-diffeomorphism from $v_* + \Omega_\gamma$ 
onto $v_* + \Omega_{\omega_e(\gamma)}$ with 
$\omega_e(\gamma) = (\gamma+r_e)/(1+2\gamma r_e + r_e^2)^{1/2}$. 
Moreover, there exists $C_\gamma \in (0,\infty)$ such that 
\bear\label{vtovprim2}
&& C^{-1}_\gamma \, |v-v_*| 
\le |\phi_e(v) - v_*|Ê\le 2 \, |v-v_*|, \\ \label{vtovprim3}
&& | \phi^{-1}_e(v') - \phi^{-1}_{e'}(v') | 
\le C_\gamma \, |e'-e|Ê\, |v'-v_*|, \\ \label{vtovprim4}
&&|J_e | \le C_\gamma, \quad |J_e^{-1}|  
\le C_\gamma, \quad  |J^{-1}_e - J^{-1}_{e'} | \le  C_\gamma \, |e'-e|
\eear
on $v_*+\Omega_\gamma$ uniformly with respect to the 
parameters $e,e' \in [0,1]$, $\sigma \in \Sph^{N-1}$, $v_* \in \R^N$. 
The same estimates hold for $\phi^*_e$. 

Finally, for any $e,e' \in [0,1]$, $\sigma \in \Sph^{N-1}$, $v_* \in \R^N$ 
and $t \in [0,1]$ there holds
\beqn\label{vtovprim5}
t \, \phi^{-1}_{e} + (1-t) \,  \phi^{-1}_{e'}  = \phi^{-1}_{e''} 
\eeqn
for some $e''$ into the segment with extremal points $e$ and $e'$. The same result holds for $\phi^*_e$. 
\end{lem}

\smallskip\noindent
{\sl Proof of Lemma \ref{v&v*tovprim}.}  
We only establish the result for the function $\phi_e$, since the proof 
for $\phi^*_e$ is similar (and even simpler).  
First, (\ref{vtovprim4}) and the fact that $\phi_e$ defines 
a $C^\infty$-diffeomorphism from $v_* + \Omega_\gamma$ 
onto $v_* + \Omega_{\omega_e(\gamma)}$ come straightforwardly from Lemma~\ref{utow} and its proof. 

Second, for any $z \in D$, $|\Phi_z(u)| = | u + |u|Ê\, z|Ê\le 2 \, |u|$ and 
$$
|\Phi_z(u)|^2 \ge |u|^2 + 2 \, |u| \, z \cdot u  + |u|^2 \, |z|^2 \ge |u|^2 \, (1-\gamma^2)
$$
for any $u \in \R^N$, $\hat u \cdot \hat z \ge \gamma$. That proves (\ref{vtovprim2}). 

Third, using the notation of Lemma~\ref{utow} we write 
$\Phi^{-1}_{r \, \sigma}(w) = (\varphi^{-1}_{w_2,r}(w_1),w_2)$
for any $w = (w_1,w_2)$, $w_1 \in \R$, $w_2 \in \R^N$, $w_2 \cdot \sigma = 0$. 
The map $(u_1,r) \mapsto \varphi_{w_2,r}(u_1)$ is smooth and 
has positive partial derivatives on $\R \times [0,1]$ 
if $w_2 \not = 0$ and on $(0,\infty) \times [0,1]$ if $w_2 = 0$.
On the one hand, we deduce that $(w_1,r) \mapsto \varphi^{-1}_{w_2,r}(w_1)$ 
is smooth and increasing in both variables and that the same holds 
for 
  \[ (w_1,e) \mapsto \frac{4}{3-e} \, \varphi^{-1}_{w_2,r_e}(w_1).  \]
The intermediate values Theorem then implies that for any 
$e \le e' \in [0,1]$, $t \in [0,1]$ there holds 
$$
t \, {4 \over 3-e} \, \varphi^{-1}_{w_2,r_e}(w_1) 
+ (1-t) \, {4 \over 3-e'} \, \varphi^{-1}_{w_2,r_{e'}}(w_1) 
= {4 \over 3-e''} \, \varphi^{-1}_{w_2,r_{e''}}(w_1) 
$$
for some $e'' \in [e,e']$ from which (\ref{vtovprim5}) follows. 

On the other hand, $r \mapsto \Phi^{-1}_{r \, \sigma}(\hat w)$ is smooth for 
any $\hat w \in \Sph^{N-1} \backslash \{-\sigma\}$ and therefore 
there exists $C_\gamma$ such that 
$|  \Phi^{-1}_{r \, \sigma} (\hat w) - \Phi^{-1}_{r \, \sigma} (\hat w) | \le C_\gamma \, |r'-r|$
uniformly for any  $\hat w \in \Sph^{N-1}$, $\hat w \cdot \sigma \ge \gamma$. 
Thanks to the homogeneity property $\Phi^{-1}_z(\lambda \, w) = \lambda \, \Phi^{-1}_z(w)$ we deduce  
$$
|  \Phi^{-1}_{r \, \sigma} (w) - \Phi^{-1}_{r \, \sigma} (w) | = |w| 
|  \Phi^{-1}_{r \, \sigma} (\hat w) - \Phi^{-1}_{r \, \sigma} (\hat w) | \le C_\gamma \, |r'-r| \, |w|,
$$
from which (\ref{vtovprim3}) follows. \qed

\smallskip \noindent
{\sl Proof of Proposition~\ref{uniqNL}.}
Let us denote $Q_f$ (resp. $Q_g$) the collision operator with collision 
rate associated with $\EE=\EE(f)$ (resp. $\EE=\EE(g)$), 
$D:=f-g$, $S:=f+g$ and $k := (1+|v|^2) \, \mbox{sgn}(D)$. The evolution equation on $D$ writes
\[
\frac{\partial}{\partial t} D = {1 \over 2}  \left[ Q_f (D,S) + Q_f (S, D) \right]
         + \left[ Q_f(g,g) - Q_g(g,g) \right] 
\]
and thus the time derivative of the $L^1 _2$ norm of $D$ is
\bean
\frac{d}{dt} \| D  \|_{L^1_2} 
&= & \frac12 \, \int_{\RR^N \times \RR^N \times \Sph^{N-1}} S D_*
         \Big[ k\big(v'_{e(f)}\big)+ k\big(v'_{*,e(f)}\big) - k - k_* \Big] 
        \, |u| \, \tilde b_{\EE(f)} \, dv \, dv_*  \, d\sigma \\
&+& \int_{\RR^N \times \RR^N \times \Sph^{N-1}} g g_* \, \Big[ k\big(v'_{e(g)}\big) - k \Big] \, 
      |u| \, \Big[ \tilde b_{\EE(f)}  - \tilde b_{\EE(g)} \Big] \, dv  \, dv_* \, d\sigma \\
&+& \int_{\RR^N \times \RR^N \times \Sph^{N-1}} g g_* \, \Big[ k\big(v'_{e(f)}\big) - k\big(v'_{e(g)}\big) \Big] \, 
        |u| \, \tilde b^+_{\EE(f)} \,  dv \, dv_* \, d\sigma \\
&+& \int_{\RR^N \times \RR^N \times \Sph^{N-1}} g g_* \, \Big[ k\big(v'_{e(f)}\big) - k\big(v'_{e(g)}\big) \Big] \, 
        |u| \, \tilde b^-_{\EE(f)}  \, dv \, dv_* \, d\sigma \\
&=:& I_1 + I_2 + I_3 + I_4,
\eean
the subscripts recalling that the post-collisional velocities 
$v'_{e(f)}$, $v'_{*,e(f)}$ and $v'_{e(g)}$ defined by (\ref{vprimvisco}) 
depend on the choice of the normal restitution coefficient $e$ and 
thus on the kinetic energies $\EE(f)$ and $\EE(g)$. 
Here we have set $\tilde b_\EE(x) = \tilde b(\EE,x) $ 
and $\tilde b^\pm_\EE(x) = \tilde b (\EE,x) \, {\bf 1}_{\pm x \ge 0}$ 
and for the sake of brevity we just write $e(h)$ instead of $e(\EE(h))$ for any function $h \in L^1_2$. 

\medskip
The first term is easily dealt with by the same arguments as in the
non-coupled case:
\[
      I_1 \le \int_{\RR^N \times \RR^N \times \Sph^{N-1}} S \, |D_*| \, (1+|v|^2)
      \, |u| \, \tilde b_{\EE(f)}  \, dv \, dv_* \, d\sigma
      \le \alpha(\EE(f)) \, \|S\|_{L^1 _3} \, \|f - g\|_{L^1 _1}.
\]

Using $|u| \, (|k| + |k(v'_{e(g)})|) \le 2 \, (1+|v|^2)^{3/2} \, (1+|v_*|^2)^{3/2}$, 
the second term $I_2$ is controlled by
\[ 
I_2 \le 2 \, \| \tilde b_{\EE(f)} - \tilde b_{\EE(g)} \|_{L^1(S^{N-1})}  \, \|g\|_{L^1_3}^2.
\]
Using now the locally Lipschitz assumption (\ref{H1-1}) 
and the fact that $\EE(f), \EE(g) \in K$ we get for some 
constant $C_K$ depending on $\tilde b$ and $K$:
\[ 
I_2 \le C_K \, |\EE(f) - \EE(g)| \, \|g\|_{L^1 _3}^2 \le C_K \, \|f - g\|_{L^1 _2} \, \|g\|_{L^1 _3}^2.
\]

As for the third term $I_3$, we use twice the change of 
variable $v \mapsto v' = \phi_e(v)$ with $v_*,\sigma$ fixed and $e = e(f)$ or $e = e(g)$. We get
\bean
I_3 
&=& \int_{\RR^N \times \Sph^{N-1}} \int_{\OO_{e(f)}} g_* \, k'  \, G(\phi^{-1}_{e(f)}) \, J^{-1}_{e(f)} \, 
       \tilde b^+_{\EE(f)}\, dv' \, dv_* \, d\sigma \\
&-& \int_{\RR^N \times \Sph^{N-1}} \int_{\OO_{e(g)}} g_* \, k'  \,  G(\phi^{-1}_{e(g)}) \, J^{-1}_{e(g)} \,
      \tilde b^+_{\EE(f)} \, dv' \, dv_* \, d\sigma,
\eean
where we have introduced the notations $G(w) := |v_*-w| \, g(w)$ 
for any $w \in \R^N$ and $\OO_e = v_* + \Omega_{\omega_e(0)}$. 
Without restriction we may assume $e(f) \le e(g)$ and 
therefore  $\OO_{e(g)} \subset \OO_{e(f)}$ since $e \mapsto \omega_e(0)$ 
is an increasing function. We then split $I_3$ as
\bean 
I_3 
&=&  \int_{\RR^N \times \Sph^{N-1}} 
\int_{\OO_{e(f)} \backslash \OO_{e(g)}} g_* \, k'  \, 
G(\phi^{-1}_{e(f)}) \, J^{-1}_{e(f)} \,  \tilde b^+_{\EE(f)} \, dv' \, dv_* \, d\sigma \\
&+& 
\int_{\RR^N \times \Sph^{N-1}} \int_{\OO_{e(g)}} g_*  \, k'  \,
\Big[  J^{-1}_{e(f)} -  J^{-1}_{e(g)} \Big] \,  G(\phi^{-1}_{e(g)}) \, 
\tilde b^+_{\EE(f)} \, dv' \, dv_* \, d\sigma  \\
&+&
\int_{\RR^N \times \Sph^{N-1}} \int_{\OO_{e(g)}} g_* \, k' \,
\Big[ G(\phi^{-1}_{e(f)}) - G(\phi^{-1}_{e(g)})  \Big] \, 
J^{-1}_{e(f)}  \,  \tilde b^+_{\EE(f)} \, dv' \, dv_* \, d\sigma \\&= & I_{3,1} + I_{3,2} + I_{3,3}.
\eean
For the first term $I_{3,1}$ we use the backward change of 
variables $v' \mapsto v = \phi^{-1}_{e(f)}(v')$ and we get 
$$
I_{3,1} =  \int_{\RR^N \times \Sph^{N-1}} \int_{\R^N} g_* \, k (v'_{e(f)}) \, G\,  \tilde b_{\EE(f)} \, 
{\bf 1}_{ 0 \le \hat u \cdot \sigma \le  \eta }  \, dv \, dv_* \, d\sigma
$$
with $\eta := \omega_{e(f)}^{-1} \circ \omega_{e(g)}(0)$. 
By inspection, the functions $(e,\gamma)\in [0,1]  \times (-1,1] 
\mapsto \omega_e(\gamma), \omega_e^{-1}(\gamma) \in (-1,1]$ 
are smooth with respect to both variables. 
From this smoothness and the fact that $ \omega^{-1}_{e} \circ \omega_e(0) = 0$ 
we deduce $| \omega^{-1}_{e'} \circ \omega_e(0)| \le C \, |e-e'|$ 
for any $e,e' \in [0,1]$ and for some constant $C \in (0,\infty)$. 
As a consequence, thanks to the Lipschitz assumption (\ref{H1-2}), we obtain
\bean
I_{3,1} &\le& \| \tilde b \|_{L^\infty}
 \int_{\RR^N \times \R^N} g \, (1+|v|)^3 \, g_* \, (1+|v_*|)^3 \left\{ \int_{\Sph^{N-1}}
{\bf 1}_{  - C \, (e(g) - e(f)) \le \hat u \cdot \sigma \le  0 } \, d\sigma \right\}  \, dv \, dv_* \\
&\le& C  \, \| g \|_{L^1_3}^2 \, |e(f) - e(g)| \le C  \, \| g \|_{L^1_3}^2 \, \| f - g \|_{L^1_2}.
\eean

For the term $I_{3,2}$, using the estimate (\ref{vtovprim4}) and the Lipschitz assumption (\ref{H1-2}), we get
        \[ 
\Big|  J^{-1}_{e(f)} -  J_{e(g)}^{-1} \Big| \le C \, |e(f) - e(g)|  \le C_K \, \|f - g \|_{L^1 _2}. 
        \]
Then doing the backward change of variable  $v' \mapsto v = \phi^{-1}_{e(g)}(v')$ and observing that $J_{e(f)}$ is bounded on $\{ u, \,\, \hat u \cdot \sigma \ge 0 \}$ thanks to (\ref{vtovprim2}), we get
     \[ 
I_{3,2} \le C_K \, \|f - g \|_{L^1 _2} \, \|g\|_{L^1 _3} ^2. 
     \]

\smallskip
We now aim to prove that for any functions $f,\,g$ which energies
$\EE_f$ and $\EE_g$ belonging to a compact $K \subset (0,\infty)$ there
exists a constant $C_K$ such that the following functional inequality
holds
\beqn\label{bddI33}
   I_{3,3} \le   C_K \, \|f - g \|_{L^1 _2} \, \|g\|_{L^1 _4} \, \|g\|_{BV_4}.
\eeqn
Let us first assume that $f$ and $g$ are smooth functions, say $f,g \in \DD(\RR^N)$.
From (\ref{vtovprim3}) and (\ref{vtovprim5}) we have
\bean
&&  \Big| G(\phi^{-1}_{e(f)}) (v')- G(\phi^{-1}_{e(g)}(v') ) \Big| \,\,  \le \\
&& \qquad \le 
 | \phi^{-1}_{e(f)}(v')  -  \phi^{-1}_{e(g)}(v') |  \int_0 ^1
 \Big|\nabla_w G ((1-t)  \phi^{-1}_{e(f)}(v')  + t \phi^{-1}_{e(g)}(v') )\Big| \, dt  \\
&& \qquad \le C \, |e(f) - e(g)| \, |v'-v|  \int_0 ^1 \Big|\nabla_w G (\phi^{-1}_{e_t}(v') )\Big| \, dt  
\eean
with $e_t \in [e(f),e(g)]$. Since then $\OO_{e(g)} \subset \OO_{e_t}$ for any $t \in [0,1]$, we deduce 
$$
I_{3,3} 
\le C \, |e(f) - e(g)|Ê\int_0^1 \int_{\RR^N \times \Sph^{N-1}} \int_{\OO_{e_t}} g_* \, | k' | \, |v'-v| \, 
 \Big|\nabla_w G (\phi^{-1}_{e_t}(v') )\Big| \, dv' \, dv \, d\sigma dt.
$$
Using finally the backward change of variable  $v' \mapsto v = \phi^{-1}_{e_t}(v')$ 
and the uniform bound (\ref{vtovprim4}) on the Jacobian $J_{e_t}$ on $v_* + \Omega_0$  we get 
$$
I_{3,3} 
\le C \, |e(f) - e(g)| \, \|g\|_{L^1 _4} \, \|g\|_{BV_4}.
$$
Therefore we obtain (\ref{bddI33}) for smooth functions. When $f,g \in
BV_4$ we argue by density, introducing two sequences of
smooth functions $(f_n)$ and $(g_n)$ which converge respectively to $f$
and $g$ in $L^1$ and are bounded in $BV_4$, we pass to the
limit $n \to \infty$ in the functionnal inequality (\ref{bddI33})
written for the functions $f_n$ and $g_n$. We then easily conclude that
(\ref{bddI33}) also holds for $f$ and $g$.  

The term $I_4$ can be dealt with similarly to the term $I_3$. Collecting all the estimates we thus get
     \[ \frac{d}{dt} \|f_t - g_t \|_{L^1 _2} \le C' _{T_*} \, \|f_t - g_t
     \|_{L^1 _2} \]
where $C' _{T_*}$ depends on $K$, $\tilde b$ and on some uniform bounds on
$\|f\|_{L^1_3}$ and $ \|g\|_{BV_4}$. This concludes the proof by a
Gronwall argument.
\qed

\smallskip

The uniqueness part of Theorem~\ref{NLcase} follows straightforwardly
from Proposition~\ref{uniqNL} and the discussion made just before
its statement.


\section{Study of the cooling process}\label{sec:CP}
\setcounter{equation}{0}
\setcounter{theo}{0}

In this section we prove the cooling asymptotic as stated in point (ii)
of Theorem \ref{Lcase}
and points (iii), (iv), (v) of Theorem \ref{NLcase}.
We first prove the collapse of the distribution function
in the sense of weak * convergence to the Dirac mass in the set
of measures.

\begin{prop}\label{ftodelta0}
Let $T_c \in (0,+\infty]$ be the time of life of the solution.
Under the sole additional assumption {\bf H2}, there holds
    \beqn\label{fwtodelta}
    f (t,.) \,\, \mathop{\wto}_{t \to T_c}  \,\, \delta_{v=0} \, \hbox{
weakly}* \mbox{ in } M^1(\R^N).
    \eeqn
\end{prop}

\smallskip\noindent
{\sl Proof of Proposition \ref{ftodelta0}}. We split the proof in two
steps.

\smallskip\noindent
{\sl Step 1.} Assume first  that $\EE  \to 0$ when $t \to T_c$.
This is always the case when $T_c <+\infty$ (since the convergence to $0$
of the kinetic energy follows from
the existence proof in this case) and it will be established
under additional assumptions on $B$ when $T_c = +\infty$
but it probably holds true under the sole assumption {\bf H2} in this case
as well.
For any $0 \le \varphi \in \DD(\RR^N \backslash \{ 0 \})$, there exists
$r > 0$ such that $\varphi = 0$ on $D(0,r)$ and then, there exists
$C_\varphi = C_\varphi(r,\| \varphi \|_\infty)$ such that $|\varphi(v)|
     \le
\, C_\varphi \, |v|^2$. As a consequence,
    $$
    \int_{\RR^N} f \, \varphi \, dv \le C_\varphi \, \EE(t) \to 0,
    $$
from which we deduce that any weak * limit $\bar \mu$ of $f$ in $M^1$
satisfies $\mbox{supp} \, \bar \mu  \subset \{ 0 \}$.
Therefore, (\ref{fwtodelta}) follows using the  conservations
(\ref{consmq}) and the energy bound (\ref{decE}).

\smallskip\noindent
{\sl Step 2.} Assume next that  $\EE \to \EE_\infty > 0$ (and thus also
$T_c = +\infty$).
Then for a fixed time $T > 0$ and for any non-negative sequence
$(t_n)$ increasing and going to $+\infty$,
there exists a subsequence $(t_{n_k})$ and a measure
$\bar \mu \in L^\infty(0,T; M^1_2)$  such that the sequence $f_k(t,v)
:= f(t_{n_k}+t,v)$ satisfies
\beqn\label{fkmu}
    f_k \wto  \bar \mu \, \hbox{ weakly}* \mbox{ in } \,
L^\infty(0,T;M^1).
\eeqn
Moreover,  for any $\varphi \in C_c(\R^N)$, there holds
    \[
    {d \over dt} \int_{\R^N} f_k \, \varphi \, dv = \langle
    Q(f_k,f_k),\varphi \rangle
    \quad\hbox{on}\quad (0,T), \]
with $\langle Q(f_k,f_k),\varphi \rangle $ bounded in $L^\infty(0,T)$.
  From Ascoli's Theorem, we get
    \[
    \int_{\R^N} f_k \, \varphi \, dv \,\, \to \,\, \int_{\R^N}  \varphi
    \, d\bar\mu(v) \quad\hbox{uniformly on}\quad [0,T].
    \]
As a consequence, for any given function $\chi_\eps \in C_c(\RR^3
\times \RR^3)$ such that
$0 \le \chi_\eps \le 1$ and $\chi_\eps(v,v_*) = 1$ for every $(v,v_*)$
such that $|v| \le \eps^{-1}$ and $|v_*| \le \eps^{-1}$ we may pass to
the limit
(using the continuity of $\Delta=\Delta(\EE,u)$ which is uniform on the
compact set
determined by $[\EE_\infty,\EE_0]$ and the support of $\chi_\eps$)
    \beqn\label{Depsn}
    \int_0^T D_\eps(f_k) \, dt \,\, \mathop{\longrightarrow}_{k \to
+\infty} \,\,
    \int_0^T \int_{\RR^N \times \RR^N} |u|^3 \Delta(\EE_\infty,u) \,
\chi_\eps(v,v_*) \, d\bar\mu \, d\bar \mu_* \, dt,
    \eeqn
where we have defined for any measure (or function) $\lambda$:
    $$
    D_\eps(\lambda) := \int_{\RR^N \times \RR^N} |u|^3 \,
    \Delta(\EE,u) \, \chi_\eps(v,v_*) \, d\lambda(v) \, d\lambda(v_*).
    $$
  From the dissipation of energy (\ref{eqdiffEE}) and
the estimate from below (\ref{hypbeta2}), there holds
    \[
    {d \over dt} \EE(t) \le - D(f) \, \, \hbox{ with } \, \,
    D(f) := \int_{\RR^N \times \RR^N} |u|^3 \,
    \Delta(\EE,u) \, f \, f_* \, dv \, dv_*,
    \]
which in turn implies  that $t \mapsto D(f(t,.)) \in L^1(0,\infty)$,
and then
    \beqn\label{Depsto0}
    \int_0^T D_\eps(f_k) \, dt \le \int_0^T D(f_k) \, dt  = \int_{t_{n_k}}^{t_{n_k}+T} D(f) \, dt
    \,\,\, \mathop{\longrightarrow}_{k \to \infty} \,\,\, 0.
    \eeqn
Gathering (\ref{Depsn}) and (\ref{Depsto0}), and letting $\eps$ goes to
$0$, we deduce that
    \[ \int_{\RR^N \times \RR^N} |u|^3 \Delta(\EE_\infty,u)  \, d\bar\mu
\, d\bar \mu_* =0
    \,\,\, \mbox{ on }\,\,\, (0,T).\]
The positivity (\ref{hypbeta1}) of $\Delta(\EE_\infty,u)$ then implies
that $\bar \mu = \bar c \, \delta_{v=\bar w}$ for some
measurable functions $\bar w: (0,T) \to \R^N$ and $\bar c: (0,T) \to
\R_+$. Moreover, from the conservation of mass and momentum
(\ref{consmq}) and the bound of energy (\ref{decE}) we deduce that
$\bar c = 1$ and $\bar w=0$ a.e. It is then classical to deduce (by the
uniqueness of the limit and the fact that it is independent on time)
that (\ref{fwtodelta}) holds.\qed

\medskip
To conclude that this weak convergence of the distribution to the
Dirac mass as time goes to infinity implies the convergence of the
kinetic
energy to $0$ ({\em i.e.}, the kinetic energy of the Dirac mass) we have
to show that no kinetic energy is escaping at infinify as $t \to T_c$.
To this purpose we put stronger
assumptions on the collision rate. The first
additional assumption {\bf H3} roughly speaking means that the energy
dissipation functional is strong enough to forbid it,
whereas the second additional assumption {\bf H4} allows to use the
uniform
propagation of moments of order strictly greater than $2$
to forbid it.

     \begin{prop}\label{EEto0}
     Let $T_c \in (0,+\infty]$ be the time of life of the solution.
     Then if either $T_c <+\infty$, or $T_c = + \infty$ and $B$
     satisfies additional assumptions {\bf H2-H3} or {\bf H2-H4}, we have
       \beqn\label{Eto0}
       \EE(t) \, \to \, 0 \quad\hbox{when}\quad t \to T_c.
       \eeqn
     \end{prop}

\smallskip\noindent
{\sl Proof of Proposition \ref{EEto0}. }
We split the proof in three steps.  \smallskip \\
{\sl Step 1. }Assume first $T_c < + \infty$. The claim follows from the
existence proof.

\smallskip\noindent
{\sl Step 2. } Assume now $T_c = + \infty$ and that $B$ satisfies
assumption {\bf H3}: (\ref{hypbeta2})-(\ref{hypbeta3}).
We argue by contradiction: assume that $\EE(t) \not\to 0$, that  is,
there exists $\EE_\infty > 0$ such that $\EE(t) \in
(\EE_\infty,\EE_{\mbox{\scriptsize{in}}})$. Reasoning as in Proposition
\ref{ftodelta0}, we
get,
for a fixed time $T > 0$ and for any sequence
$(t_n)$ increasing and going to infinity,
that there exists a subsequence $(t_{n_k})$ and a measure
$\bar \mu \in L^\infty(0,T; M^1_2)$  such that the function
$f_k(t,v) := f(t_{n_k}+t,v)$ satisfies (\ref{fkmu})
and
    \beqn\label{Depsnpsi}
    \int_0^T D_\eps ^0(f_k) \, dt \to \int_0^T D_\eps ^0(\bar\mu) \, dt,
    \eeqn
where we have defined for any measure (or function) $\lambda$:
    $$
    D_\eps ^0 (\lambda) := \int_{\RR^N \times \RR^N} |u|^3 \,
    \psi(|u|) \, \chi_\eps(v,v_*) \, d\lambda(v) \, d\lambda(v_*).
    $$
  From the dissipation of energy (\ref{eqdiffEE}) and
the estimate from below (\ref{hypbeta2}), there holds
\beqn\label{diffED0}
{d \over dt} \EE(t) \le - D^0(f) \quad\hbox{with}\quad
\quad D^0(f) := \int_{\RR^N \times \RR^N} |u|^3 \,
\psi(|u|) \, f \, f_* \, dv \, dv_*,
\eeqn
which in turn implies  that $t \mapsto D^0(f(t,.)) \in L^1(0,\infty)$,
and then
\beqn\label{Depsto0psi}
\int_0^T D_\eps ^0(f_k) \, dt \le \int_0^T D^0(f_k) \, dt = \int_{t_{n_k}}^{t_{n_k}+T}
D^0(f) \, dt
\,\,\, \mathop{\longrightarrow}_{k \to \infty} \,\,\, 0.
\eeqn
Gathering (\ref{Depsnpsi}) and (\ref{Depsto0psi}), and letting $\eps$
goes to
$0$, we deduce that  $D^0(\bar \mu) = 0$ on $(0,T)$. The
positivity of $\psi$ implies as in Proposition \ref{ftodelta0}
that $\mbox{supp} \, \bar \mu  \subset \{ 0 \}$ and $\bar\mu \delta_{v=0}$.
As this limit is unique and independent on time we deduce that
(\ref{fwtodelta}) holds.

\smallskip
Now, on the one hand, taking $R = \sqrt{\EE_\infty/2}$ there holds
    \beqn\label{BRcfv2}
    \int_{B_R^c} f \, |v|^2 \, dv = \int_{\R^N} f \, |v|^2 \, dv -
    \int_{B_R} f \, |v|^2 \, dv
    \ge \EE_\infty - R^2 \ge \EE_\infty/2
    \eeqn
for any $t \ge 0$. On the other hand, for $T$ large enough, there holds
thanks to (\ref{fwtodelta})
    \beqn\label{BRf}
    \int_{B_{R/2}} f \, dv \ge {1 \over 2} \quad\hbox{for any}\quad t \ge
T.
    \eeqn
Remarking that on $B_{R/2} \times B_R^c$ there holds, thanks to
(\ref{hypbeta3}),
    \beqn\label{u3pR}
    |u|^3 \, \psi(|u|) \ge {|v_*|^3 \over 8} \, \psi \left(
    {|v_*| \over 2} \right)
    \ge \psi_R \, {|v_*|^2 \over 4},
    \eeqn
we may put together (\ref{diffED0})-(\ref{u3pR}) and we get thanks to
(\ref{BRcfv2}) and  (\ref{BRf})
\bean
{d \over dt} \EE(t) &\le& - \int_{B_{R/2}} \! \int_{B_R^c} |v - v_*|^3
\, \psi(|v-v_*|) \, f \, f_* \, dv dv_* \\
&\le& - {\psi_R \over 4}\int_{B_{R/2}} f \, dv  \int_{B_R^c} f_* \,
|v_*|^2 \,  dv_* \le  - {\psi_R \over 4} \,  {1\over 2} \, {\EE_\infty
\over 2}
\eean
for any $t \ge T$. This implies that $\EE$ becomes negative in finite
time and we get a
contradiction.

\smallskip\noindent
{\sl Step 3. }
Finally,  assume that $T_c = + \infty$ and $B$ satisfies assumption
{\bf H4}.
On the one hand, thanks to (\ref{MomentUnif1}), there holds
    $$
    \sup_{[0,\infty)} \int_{\RR^N} f(t,v) \, |v|^3 \, dv < \infty.
    $$
On the other hand, arguing as in Step 2, we obtain (keeping the same
notations) that (\ref{fkmu})  and then (from the uniform bound in $L^1
_3$)
    $$
    \EE(f_k) \to \bar\EE=\EE(\bar\mu) \ \quad\hbox{and  }\quad
     D(\bar\mu) = 0.
    $$
The dissipation of energy vanishing implies that
    $$
    |u|^3 \, \bar \mu \, \bar \mu_* \equiv 0 \quad\hbox{or}\quad
\Delta(\bar\EE,u)
    \,\,\hbox{is not positive on }  (0,T) \times \R^{2N}.
    $$
In the first case we deduce that $\bar\mu = \delta_{v=0}$ as in Step 2
and then $\bar\EE = \EE( \delta_{v=0}) = 0 $. In the second case we
deduce, from (\ref{hypbeta1}), that $\bar\EE$ is not positive. In both
case, there exists $\tau_k$ such that $\tau_k \to \infty$ and
$\EE(\tau_k) \to 0$ and therefore (\ref{EEto0}) holds since $\EE$ is
decreasing. \qed

\medskip

\smallskip
Now we turn to some criterions for the cooling process to occur or not
in finite time.

    \begin{prop}\label{alphaj}
    Assume that $\alpha$ is bounded near $\EE = 0$,
    and $j_\EE$ converges to $0$ as $\e \to 0$ uniformly
    near $\EE = 0$, then $T_c =+\infty$.
    \end{prop}

\smallskip\noindent
{\sl Proof of Proposition \ref{alphaj}}.
It is enough to remark that, thanks to the hypothesis made on $\alpha$
and $j_\EE$, the {\it a priori} bound in Orlicz norm that one deduces
from
(\ref{dnormLLambdadt}) as in Corollary \ref{cor:Orl} extends to all
times:
     \[ \forall \, t \ge 0 \qquad  \|f_t\|_{L^\Lambda} \le
\|f_{\mbox{\scriptsize{in}}}\|_{L^\Lambda} \, \exp \left( C \, \|
f_{\mbox{\scriptsize{in}}} \|_{L^1_2} \, t
\right) \]
for some constant $C$ depending on the collision rate. It shows that
the energy cannot
vanish in finite time.
\qed

    \begin{prop}\label{deltale0}
    Assume that $B$ satisfies {\bf H4}, that for some increasing and
    positive function  $\Delta_0$ there holds  $\Delta(\EE,u) \le
    \Delta_0(\EE)$ for any $u \in \R^N$, $\EE \ge 0$, and that
$f_{\mbox{\scriptsize{{\em in}}}}\,
    e^{r \, |v|^\eta}\in L^1$ for some $r > 0$ and $\eta \in (1,2]$, then
    $T_c = +\infty$.
    \end{prop}

\smallskip\noindent
{\sl Proof of Proposition \ref{deltale0}}.
  From the dissipation of energy (\ref{eqdiffEE}),
the bound on $\Delta$ and the decay of
the energy (\ref{decE}), we have
    $$
    {d \EE \over dt} \ge -  \Delta_0(\EE_{\mbox{\scriptsize{in}}})  \,
\int_{\RR^N} \!
    \int_{\RR^N} \,
    f \, f_* \,  \, |u|^3 \, dvdv_* =: -
\Delta_0(\EE_{\mbox{\scriptsize{in}}})  \,  (I_{1,R} +
    I_{2,R})
    $$
where
    \[ \left\{
       \begin{array}{l}\displaystyle
       I_{1,R} := \int_{\RR^N \times \RR^N} |u|^3 \, {\bf 1}_{\{|u| \le
R\}} \, f \, f_* \, dv \, dv_*,
       \vspace{0.2cm} \\ \displaystyle
       I_{2,R} := \int_{\RR^N \times \RR^N} |u|^3 \, {\bf 1}_{\{|u| \ge
R\}} \, f \, f_* \, dv \, dv_*.
       \end{array}
       \right. \]
On the one hand, for any $R > 0$, we have using (\ref{consmq})
      \[
      I_{1,R} \le R \int_{\RR^N \times \RR^N} |u|^2 \,  f \, f_*
       \, dv \, dv_*=  2 \, R \, \EE. \]
On the other hand, we infer from Proposition~\ref{Y3uniform} (since
$B$ satisfies {\bf H4}) that
    $$
    \sup_{t \in [0,T_c)} \int_{\R^N} f(t,v) \, e^{2 \, r' \, |v|^\eta} \,
dv
    \le C_1
    $$
for some $r',C_1 \in (0,\infty)$. Therefore
    \bean
    I_{2,R} &\le& \int_{\RR^N \times \RR^N} (4 \,|v|^3 + 4 \, |v_*|^3) \,
2 \,
    {\bf 1}_{\{|v| > R/2\}} \, f \, f_* \, dv \, dv_* \\
    &\le& {8 \, e^{-r' \, R^\eta}} \int_{\RR^N} (1 + |v|^3) \,
    e^{r' \, |v|^\eta} \, f \, dv
    \int_{\RR^N} (1 + |v_*|^3) \, f _* \, dv_*  \le C_2 \, e^{-r'
\,R^\eta}.
    \eean
Gathering these three estimates, we deduce
    $$
    {d \over dt} \EE \ge - C_3 \, R \, \EE - C_3 \, e^{- r' \, R^\eta},
    $$
which in turns implies, thanks to a Gronwall argument,
    $$
    \forall \, R > 0, \quad \inf_{t \in [0,T]} \EE(t) \ge
\EE_{\mbox{\scriptsize{in}}} \,
    e^{- C_3 \, R \, T} - {e^{-r' \, R^\eta} \over R}.
    $$
We conclude that $\EE(t) > 0$ for any $t \in [0,T]$ and any fixed $T >
0$,  choosing $R$ large enough (using that $\eta >1$).
\qed

\begin{prop}\label{deltage-1/2}
Assume $\Delta(\EE,u) \ge \Delta_0
\, \EE^\delta$ with $\Delta_0 > 0$  and $\delta < -1/2$, then $T_c <
+\infty$.
\end{prop}

\smallskip\noindent
{\sl Proof of Proposition \ref{deltage-1/2}}. On the one hand, from the
dissipation of energy (\ref{eqdiffEE}) and the bound on $\Delta$,  we
have
$$
{d \EE \over dt} \le -  \Delta_0 \, \EE^\delta  \, \int_{\RR^N} \!
\int_{\RR^N} \,
f \, f_* \,  \, |u|^3 \, dv\, dv_*.
$$
On the other hand,  from Jensen's inequality and the conservation of
mass
and momentum, there holds
$$
\int_{\RR^N} \! \int_{\RR^N} \,
f \, f_* \,  \, |u|^3 \, dvdv_*\ge \left(\int_{\RR^N} \! \int_{\RR^N} \,
f \, f_* \,  \, |u|^2 \, dvdv_*\right)^{3/2} = (2 \, \EE)^{3/2}.
$$
Gathering these two estimates, we get
      $$
      {d \over dt} \EE \le - \Delta_0 \, \EE^{\delta+3/2}
      $$
and $\EE$ vanishes in finite time. \qed

\bigskip
\appendix

\section*{Appendix: Some facts about Orlicz spaces}
\addcontentsline{toc}{section}{Appendix: Some facts about Orlicz spaces}

\def\theequation {{A.\arabic{equation}}}
\def\thetheo {{A.\arabic{theo}}}
\def\thesection{}
\setcounter{equation}{0}
\setcounter{theo}{0}

The goal of this appendix is to gather some results
about Orlicz spaces in order to make this paper
as self-contained as possible. The definition and H\"older's
inequality are recalls of results which can be found
in~\cite{RR91} for instance.
We also state and prove a simple formula for the differential of
Orlicz norms,
which is most probably not new, but for which we were not able
to find a reference.

\subsection*{Definition}

We recall here the definition of Orlicz spaces on $\R^N$ according to
the Lebesgue measure.
Let $\Lambda: \RR_+ \to \RR_+$ be a function $C^2$ strictly increasing,
convex, such that
      \beqn \label{eq:hypLamb1}
      \Lambda(0) = \Lambda'(0) = 0,
      \eeqn
      \beqn \label{eq:hypLamb2}
      \forall \, t \ge 0, \ \ \ \Lambda(2 \, t) \le c_\Lambda \,
\Lambda(t),
      \eeqn
for some constant $c_\Lambda >0$, and which is superlinear, in the
sense that
      \beqn \label{eq:hypLamb3}
      {\Lambda (t) \over t} \mathop{\longrightarrow}_{t \to + \infty} +
\infty.
      \eeqn
We define $L^\Lambda$ the set of measurable functions $f : \R^N
\rightarrow \R$ such that
      \[ \int_{\RR^N} \Lambda \bigl( |f(v)| \bigr) \, dv < + \infty. \]
Then $L^\Lambda$ is a Banach space for the norm
      \[ \|f\|_{L^\Lambda} = \inf \left\{ \lambda > 0 \ | \ \ \int_{\R^N}
\Lambda \left( {|f(v)| \over
         \lambda}  \right) \, dv \le 1 \right\} \]
and it is called the {\it Orlicz space} associated with $\Lambda$.
The proof of this last point can be found in~\cite[Chapter~III,
Theorem~3]{RR91}.
Note that the usual Lebesgue spaces $L^p$ for $1 \le p < +\infty$
are recovered as particular cases of this definition for
$\Lambda (t) = t^p/p$.

Let us mention that for any $f \in L^1(\R^N)$,
a refined version of the De la Vall\'ee-Poussin Theorem
\cite[Proposition~I.1.1]{Le} (see also \cite{LMb,LM02}) guarantees that
there exists a function $\Lambda$ satisfying all the properties above
and such that
      \[ \int_{\R^N} \Lambda(|f(v)|) \, dv < +\infty.\]

\subsection*{H\"older's inequality in Orlicz spaces}

Let $\Lambda$ be a function $C^2$ strictly increasing, convex satisfying
the assumptions~(\ref{eq:hypLamb1}), (\ref{eq:hypLamb2}) and
(\ref{eq:hypLamb3}), and $\Lambda^*$ its {\it complementary Young
function},
given (when $\Lambda$ is $C^1$) by
      \[ \forall \, y \ge 0, \ \ \ \Lambda^*(y)=y (\Lambda')^{-1}(y) -
\Lambda((\Lambda')^{-1}(y)). \]
It is straightforward to check that $\Lambda^*$ satisfies the same
assumptions as $\Lambda$. Recall Young's inequality
      \beqn\label{YoungIneg}
      \forall \ x,\, y \ge 0, \qquad x \, y \le \Lambda(x) +
\Lambda^*(y).
      \eeqn
Then one can define the following norm on the Orlicz space
$L^{\Lambda^*}$:
      \[ N ^{\Lambda^*} (f) = \sup \left\{ \int_{\RR^N} |fg| \, dv \, ;
           \ \int_{\RR^N} \Lambda (|g|) \, dv \le 1 \right\}. \]
One can extract from~\cite[Chapter~III, Section~3.4, Propositions~6
and~9]{RR91} the following
result
     \begin{theo}\label{theo:holderOrl}
     (i) We have the following H\"older's inequality for any $f \in
L^\Lambda$, $g \in L^{\Lambda^*}$:
       \beqn \label{eq:holdOrl}
       \int_{\RR^N} |fg| \, dv \le \|f\|_{L^\Lambda} \, N ^{\Lambda^*}
(g).
       \eeqn
      (ii) There is equality in~(\ref{eq:holdOrl}) if and only if there
is
      a constant $0<k^*<+\infty$ such that
        \begin{equation} \label{eq:holdegalOrl}
        \left( \frac{|f|}{\|f\|_{L^\Lambda}} \right) \left( \frac{k^*
        |g|}{N ^{\Lambda^*} (g)} \right)
        = \Lambda \left(\frac{|f|}{\|f\|_{L^\Lambda}}\right)
          + \Lambda^* \left(\frac{k^* |g|}{N ^{\Lambda^*} (g)} \right)
        \end{equation}
      for almost every $v \in \R^N$.
      \end{theo}

\subsection*{Differential of Orlicz norms}

In order to propagate bounds on Orlicz norms along the
flow of the Boltzmann equation, we shall need a formula
for the time derivative of the Orlicz norm.

     \begin{theo}\label{theo:diffOrlNorm}
     Let $\Lambda$ be a function $C^2$
     strictly increasing, convex satisfying~(\ref{eq:hypLamb1}),
     (\ref{eq:hypLamb2}), (\ref{eq:hypLamb3}), and let $0 \le f \in C^1([0,T],L^\Lambda)$ 
     such that $f(t,\cdot) \not\equiv 0$ for all $t \in [0,T]$.
     Then we have
       \beqn
       \frac{d}{dt} \|f_t\|_{L^\Lambda} 
       = \left[ N ^{\Lambda^*} \left( \Lambda' \left(
       \frac{|f|}{\|f\|_{L^\Lambda}} \right) \right) \right]^{-1}
       \int_{\R^N} \partial_t f \, \Lambda'\left(
       \frac{|f|}{\|f\|_{L^\Lambda}} \right) \, dv.
       \eeqn
     \end{theo}
{\sl Proof of Theorem~\ref{theo:diffOrlNorm}.}
From~\cite[Chapter~III, Proposition~6]{RR91}), our assumptions on
$\Lambda$ imply that
      \beqn
      \int_{\R^N} \Lambda \left( \frac{|f|}{\|f\|_{L^\Lambda}} \right) \,
dv = 1
      \eeqn
for all $0 \not= f \in L^\Lambda$.
By differentiating this quantity along $t$ we deduce:
      \[ 0 = \int_{\R^N} \partial_t f \, \Lambda'\left(
        \frac{|f|}{\|f\|_{L^\Lambda}} \right) \, dv
              - \frac{1}{\|f_t\|_{L^\Lambda}} \frac{d}{dt}
            \|f_t\|_{L^\Lambda} \,
               \int_{\R^N} f \, \Lambda'\left(
\frac{|f|}{\|f\|_{L^\Lambda}}  \right) \, dv.  \]
Now using the case of equality in H\"older's
inequality~(\ref{eq:holdOrl}) we have
      \[ \int_{\R^N} f \, \Lambda'\left( \frac{|f|}{\|f\|_{L^\Lambda}}
\right) \,
dv
         = \|f\|_{L^\Lambda} \, N ^{\Lambda^*} \left( \Lambda' \left(
         \frac{|f|}{\|f\|_{L^\Lambda}} \right) \right) \]
since the equality~(\ref{eq:holdegalOrl}) is trivially satisfied with
      \[ g = \Lambda' \left( \frac{|f|}{\|f\|_{L^\Lambda}} \right) \]
and $k^* = N ^{\Lambda^*} (g)$, using that
      \[ x y = \Lambda(x) + \Lambda^*(y) \]
as soon as $y = \Lambda'(x)$. This concludes the proof. \qed

\bigskip
\medskip

\noindent
{\bf{Acknowledgment}}: The authors thank F. Filbet, P. Lauren\c cot
and V. Panferov for fruitful remarks and discussions.
Support by the European network HYKE, funded by
the EC as contract HPRN-CT-2002-00282, is acknowledged.

\end{document}